\documentclass[12pt,twoside,reqno]{amsart}
\usepackage[colorlinks=true,citecolor=blue]{hyperref}
\usepackage{mathptmx, amsmath, amssymb, amsfonts, amsthm, mathptmx, enumerate, color}
\usepackage{graphicx}
\usepackage{epstopdf}

\usepackage{bm}

\newtheorem{theorem}{Theorem}[section]

\theoremstyle{definition}
\newtheorem{definition}{Definition}[section]

\numberwithin{equation}{section}

\newcommand{\jbe}{\begin{equation}}
\newcommand{\jen}{\end{equation}}

\newcommand{\kg}{\mbox{}\hspace{0.3in}}

\newtheorem{cor}{Corollary}[section]
\newtheorem{ex}{Example}[section]

\newtheorem{prop}{Proposition}[section]

\begin{document}


\bigskip
\baselineskip 16pt
{\footnotesize
\begin{center}
{\Large\bf Exact Penalty Algorithm of Strong Convertible Nonconvex Optimization }

\bigskip Min Jiang$^{1}$, Rui Shen$^{2}$,Zhiqing Meng$^{1}$, Chuangyin Dang$^3$\\

1. School of Management, Zhejiang University of Technology, Hangzhou, Zhejiang, 310023, China\\

2. School of Economics, Zhejiang University of Technology, Hangzhou, Zhejiang, 310023, China\\

3. Department of System Engineering and Engineering Management,City University of Hong Kong, Kowloon, Hong Kong
\end{center}

{\bf Abstract}
This paper defines a strong convertible nonconvex(SCN) function for solving the unconstrained optimization problems with the nonconvex or nonsmooth(nondifferentiable) function. First, many examples of SCN function are given, where the SCN functions are nonconvex or nonsmooth. Second, the operational properties of the SCN functions are proved, including addition, multiplication, compound operations and so on. Third, the SCN forms of some special functions common in machine learning and engineering applications are presented respectively where these SCN function optimization problems can be transformed into minmax problems with a convex and concave objective function. Fourth,a minmax optimization problem of SCN function and its penalty function are defined. The optimization condition,exactness and stability of the minmax optimization problem are proved. Finally, an algorithm of penalty function to solve the minmax optimization problem and its convergence are given. This paper provides an efficient technique for solving unconstrained nonconvex or nonsmooth(nondifferentiable) optimization problems to avoid using subdifferentiation.

{\bf  Key words:} Nonconvex nonsmooth optimization; strong convertible nonconvex function; SCN function; Exact penalty function.
}

\section{Introduction}
In this paper the following unconstrained optimization (convertible nonconvex optimization, CNO)  with a strong  convertible nonconvex function(SCN) is considered:
\begin{eqnarray*}
\mbox{(CNO)}\qquad & \min\; & f(\bm x) \\
& \mbox{s.t.}\; & \bm x\in R^n,
\end{eqnarray*}
where $f: R^n \rightarrow R$ is neither convex nor smooth. In machine learning, there are many nonconvex, nonsmooth, non-Lipschitz and discontinuous optimization problems in \cite{Li1,Mohri,Wang,Zhang}. So, to solve these problems, theoretical tools of nonsmooth and nonconvex functions are needed, such as the subdifferentiable, general convex, smoothing and so on in \cite{Bagirov1,Chieu,Clarke,Grant,Rockafellar}. A new nonconvex function is defined in this paper, which is called the SCN function in Definition 2.1, where  the SCN function is a nonconvex nonsmooth function form that can be transformed into a convex smooth function with convex equality constraints. The SCN function somewhat relates to upper -$UC^k$ function \cite{Bougeard,Daniilidis,Hare2,Ngai,Rockafellar,Spingarn} and factorable nonconvex function \cite{Bunin,Granvilliers,He,Jackson,Lundell1,McCormick,Tawarmalani2,Serrano}.

The lower(upper)-$C^k$ function was suggested by Professor R. T. Rockafellar\cite{Rockafellar}. The class of lower-$C^1$ functions is first introduced by Spingarn in \cite{Spingarn}.  In his work, Spingarn
showed that these functions are (Mifflin) semi-smooth and Clarke regular and are characterized by a generalized monotonicity property of their subgradients, called submonotonicity. The definition of lower(upper)-$C^k$ function is given as follows\cite{Daniilidis}.
\begin{definition}
Let $U$ be an open subset of $R^n$ and $k\in N$.
 Function $f: U\to R^1$ is called lower-$C^k$(for short, $LC^k$), if for every $x_0\in U$ there exist $\delta> 0$, compact topological space $S$, and a jointly continuous function $F: B(x_0, \delta)\times S\to R^1$ satisfying
$$f(x) = \max\limits_{s\in S} F(x, s),\mbox{    for all  }  x\in B(x_0, \delta),$$
such that all derivatives of $F$ up to order $k$ with respect to $x$ exist and are jointly continuous. If $-f$ is lower-$C^k$, then $f$ is called upper-$C^k$ function.
\end{definition}

The lower(upper)-$C^k$ function is nonconvex or nondifferentiable, but it is locally Lipschitz approximately convex functions in \cite{Daniilidis}.
Research on the lower(upper)-$C^k$ functions is done on subdifferentiation and optimization in \cite{Hare1,Hare2,Hare3,Hare4}.  The Moreau envelopes  $erf$:
$$erf(x) := \inf\limits_{w}\{ f(w)+\frac{r}{2}|w-x|^2\} $$
is lower-$C^2$ in \cite{Bougeard,Hare2,Ngai} such that subdifferential of the lower(upper)-$C^k$ functions can solve nonconvex optimization by prox-regularity and the proximal mapping(operator) in \cite{Hare3}.  Chieu et al. proved
second-order necessary and sufficient conditions for lower-$C^2$ functions to be convex and strongly convex in \cite{Chieu}.

Some methods for non-smooth non-convex optimization programs with lower(upper)-$C^k$ functions have been studied in \cite{Dao,Hare5,Hare6,Noll}. Dao developed a nonconvex bundle method based on the downshift mechanism and a proximity control management technique to solve nonconvex nonsmooth constrained optimization problems, where he  proved its global convergence in the sense of subsequences for both classes of lower-$C^1$ and upper-$C^1$ in \cite{Dao}.
 Hare et al.  studied two  proximal bundle methods for nonsmooth nonconvex optimization in \cite{Hare5,Hare6} by proximal mapping on lower-$C^2$ functions. Noll defined a first-order model of $f$ as an extend case of lower-$C^k$ function and presented a bundle method in \cite{Noll} as follows.
 \begin{definition}
 A function $\phi: R^n\times R^n\to R^1$ is called a first-order model of $f$ on $\omega\subset
R^n$, if $\phi(\cdot, x)$ is convex for every fixed $x\in \omega$, and if the following axioms are satisfied:

(M1) $\phi(x, x) = f(x)$ and $\partial_1\phi(x, x)\subset \partial f(x)$.

(M2) For every sequence $y_j\to x$ there exists $\epsilon_j\to 0^+$ such that $f(y_j)\leq  \phi(y_j, x)+\epsilon_j\|y_j-x\|$ for all $j\in {N}$.

(M3) For sequences $y_j\to y \in R^n$ and $x_j\to x$ in $\omega$ one has $f(y_j)\leq \limsup_{j\to \infty}\leq  \phi(y, x)$ for all $j\in {N}$.
\end{definition}
Clearly, if $f$ a first-order model, $f$ is not necessarily lower-$C^k$, and the reverse is not necessarily true.

On the other hand, the branch-and-bound method in conjunction with underestimating convex problems had been proved as an effective method to solve global nonconvex optimization problems in \cite{Al-Khayyaltt,Bao1,Tawarmalani1}. Almost all the methods used to solve nonconvex optimization are to construct many convex relaxation subproblems with convex envelopes and convex underestimating, as in
\cite{Bao1,Sherali2,Serrano,Tawarmalani2}. Based on this idea, the factorable programming technique, one of the most popular approaches for constructing convex relaxations of nonconvex optimization problems including problems with convex-transformable functions, was given in \cite{McCormick}. Due to its simplicity, factorable programming technique is included in most global optimization packages such as BARON(1996), ANTIGONE(2014), etc\cite{Nohra}. But,  Nohra and Sahinidis(2018) pointed out that a main drawback of factorable programming technique is that it often results in large relaxation gaps in \cite{Nohra}.

In 1976, McCormick(1976)\cite{McCormick} first defined factorable nonconvex function, but factorable nonconvex function is not necessarily lower-$C^1$, such as $f(x)=|x|^{0.1}+|x+1|^{0.2}$ on $x\in R^1$, because $f(x)=|x|^{0.1}+|x+1|^{0.2}$ is not locally Lipschitz in \cite{Chen1}.
In fact, the factorable nonconvex functions in \cite{Jackson,Lundell1,McCormick,Tawarmalani2} may be special CN functions (see Definition 2.2). In recent years,  research on nonconvex factorable programming further shows its effectiveness in solving the global optimization, as shown in \cite{Bunin,Granvilliers,He,Serrano}.

In order to solve (CNO), it is  meaningful to set up a new theory and an algorithm for (CNO). A new nonconvex function is defined in this paper, which is called the strong convertible nonconvex(SCN) function. The SCN function is a nonconvex or nonsmooth(nondifferentiable) function form that can be transformed into a convex-concave smooth objective function with convex function constraints, as shown in Definition 2.1. That is, by Definition 2.1, there are functions $g,g_i(i=1,2,\cdots,s)$ and $h_j(j=1,2,\cdots,r)$ such that
$$\min\limits_{\bm x} f(\bm x)=\min\limits_{(\bm x,\bm y)}\max\limits_{\bm z} \{g(\bm x,\bm y,\bm z)\mid (\bm x,\bm y,\bm z)\in X_c(f)\}.$$
 Hence, (CNO) can be transformed into a minmax problem. We can find  a large number of examples where nonconvex or nonsmooth functions are SCN functions. Because the minmax problem is convex-concave and differentiable, the SCN function makes it easier to solve the global optimal solution to (CNO).

Since last year, we have being studied a transformation technology of nonconvex and nonsmooth functions, called convertible nonconvex(CN) functions\cite{Jiang}. The CN function is a nonconvex or nonsmooth function form that can be transformed into a convex  smooth objective function with convex function constraints in Jiang(2021).
That is, there are convex functions $g,g_i,h_j:R^n\times R^{m}(i=1,2,\cdots,s,j=1,2,\cdots,r)$ such that
$$\min\limits_{\bm x} f(\bm x)=\min\limits_{(\bm x,\bm y)}\{g(\bm x,\bm y)\mid g_i(\bm x,\bm y)= 0,i=1,2,\cdots,s,h_j(\bm x,\bm y)= 0,j=1,2,\cdots,r\}.$$
In order to solve (CNO), Jiang et al(2021) have discussed optimal condition, Lagrangian dual and an algorithm for the unconstrained convertible nonconvex optimization problems under two different CN function forms respectively.
 Different from Jiang(2021)'s paper, we define a SCN function in this paper. The SCN function may be a CN function, but the reverse is not necessarily true. The research  hereinafter shows that  SCN form of  SCN function transformation is different from CN form of CN function transformation.

The main contribution of this paper is as follows:  methods of transforming nonconvex nonsmooth functions into convex smooth objective functions with convex constraint are proposed, and a nonconvex optimization problem can be transformed into a minmax problem, where its objective function is a convex-concave function and  its constraint is a convex set. The major advantage of SCN function is that it overcomes the disadvantage of smoothing, while ensuring that the objective function of the transformed minmax problem is convex and concave and the constraint set is  convex. Because there are many good algorithms to solve such minmax problems, the SCN function technique will become an effective method to solve nonconvex and nonsmooth functions. Different from all other current nonconvex and nonsmooth relaxation techniques, the forms of SCN function and CN function are a set of convex and smooth functions. The (CNO) composed of SCN (CN) function is equivalent to a smooth constrained optimization problem. However, the main drawback of SCN is that many variables are added into the transformed SCN form. This paper defines penalty function of a minmax problem of a SCN function. The penalty function is exact and stability under some condition.

The remainder of the paper is organized as follows. In Section 2,  a SCN function is defined. Some examples are given. In section 3, the operational properties of the SCN function are proved, including addition, multiplication and compound operations. In section 4, the SCN forms of some specially structured functions are discussed. In section 5, a minmax optimization problem of a SCN function is defined. The equivalence of optimality condition, exactness and stability of this minmax problem are proved. In section 6, the conclusion is given.

\section{Strong Convertible Nonconvex Function}

In this section,  SCN function is defined. Some examples are given to show that nonconvex or discontinuous functions are differentiable SCN functions.
\begin{definition}
Let $S_1\subset R^n, S_2\subset R^{m_1}$ and $S_3\subset R^{m_2}$ be convex sets and $S=S_1\times S_2\times S_3$. Let functions $g_i,h_j:R^n\times R^{m_1}\times R^{m_2}\to R^1$ ($m_1\geq 0,m_2\geq 0$,$i=1,2,\cdots,s$ and $j=1,2,\cdots,r$) be convex and differentiable on $S$, where $g_i$ is nonlinear or linear at $i=1,2,\cdots,s$ and $h_j$ is linear at $j=1,2,\cdots,r$ on $S$ respectively. Let function $g: S\to R^1$ be  convex and differentiable on $S_1\times S_2$ and be  concave and differentiable on $S_3$. Let function $f: R^n \rightarrow R^1$ be nonconvex or nonsmooth.
Let  set
\begin{eqnarray}
  X_c(f)=\{(\bm x,\bm y,\bm z)\in S&\mid & g_i(\bm x,\bm y,\bm z)\leq 0,i=1,2,\cdots,s;\nonumber\\
   &&h_j(\bm x,\bm y,\bm z)=0, j=1,2,\cdots,r,\}. \label{d1}
\end{eqnarray}
 If  there is some $(\bm x,\bm y,\bm z)\in X_c(f)$ such that
\begin{eqnarray}
 f(\bm x)=g(\bm x,\bm y,\bm z),  \label{d2}
\end{eqnarray}
and
\begin{eqnarray}
 f(\bm x)=\min\limits_{\bm y} \max\limits_{\bm z}\{ g(\bm x,\bm y,\bm z) \mid (\bm x,\bm y,\bm z)\in X_c(f)\}. \label{d3}
\end{eqnarray}
hold, then $f$ is called a strong convertible nonconvex(SCN) function on $S$ (when $S=R^n\times R^{m_1}\times R^{m_2}$, the term "on $S$" is omitted).
$[g:g_1,g_2,\cdots,g_s;h_1,h_2,\cdots,h_r]$ is called a strong convertible nonconvex(SCN) form of $f$ on $S$, briefing as $f=[g:g_1,g_2,\cdots,g_s;h_1,h_2,\cdots,h_r]$. $X(f)$ is called a SCN set. For $f$, the number of its SCN form is more than one, and the number of its SCN set is more than one.

If $-f$ is  a strong  convertible nonconvex(SCN) function on $S$, then $f$ is called a strong converse convertible nonconvex(SCCN) function on $S$, $[g:g_1,g_2,\cdots,g_s;h_1,h_2,\cdots,h_r]$ is a SCN form of $-f$.

Let set
\begin{eqnarray}
 X(f)=\{(\bm x,\bm y,\bm z)|(\bm y,\bm z)=\arg\min\limits_{\bm y} \max\limits_{\bm z}\{ g(\bm x,\bm y,\bm z) \mid (\bm x,\bm y,\bm z)\in X_c(f)\}\}. \label{d4}
\end{eqnarray}

\end{definition}

It is clear that $X(f)\subset X_c(f)$ and $ X_c(f)$ is convex set.

When $f$ is  a SCCN function on $S$, \eqref{d3} becomes
\begin{eqnarray*}
 -f(\bm x)=\min\limits_{\bm y} \max\limits_{\bm z}\{ g(\bm x,\bm y,\bm z) \mid (\bm x,\bm y,\bm z)\in X_c(f)\}, \label{d5}
\end{eqnarray*}
i.e.
\begin{eqnarray*}
 f(\bm x)=\max\limits_{\bm y} \min\limits_{\bm z}\{ -g(\bm x,\bm y,\bm z) \mid (\bm x,\bm y,\bm z)\in X_c(f)\},
\end{eqnarray*}
where $g$ is convex and differentiable on $S_1\times S_2$ and is concave and differentiable on $S_3$. So, we have
\begin{eqnarray*}
\max\limits_{\bm x\in S_1} f(\bm x)=\max\limits_{\bm x}\max\limits_{\bm y} \min\limits_{\bm z}\{ -g(\bm x,\bm y,\bm z) \mid (\bm x,\bm y,\bm z)\in X_c(f)\},
\end{eqnarray*}

In Definition 2.1, it is possible  there is a special cases of the SCN forms of $f$.


 There is a SCN function where  $\bm z$ in the SCN form of $f$ disappears and  $m_2=0$. So,
$g(\bm x,\bm y,\bm z)=g(\bm x,\bm y)$ and
\begin{eqnarray*}
 f(\bm x)= \min\limits_{\bm y} g(\bm x,\bm y)\ \ \mbox{s.t.}\ (\bm x,\bm y)\in X_c(f).
\end{eqnarray*}
That is $f(\bm x)=[g(\bm x,\bm y):g_1(\bm x,\bm y),g_2(\bm x,\bm y),\cdots,g_s(\bm x,\bm y)]$. For all $(\bm x,\bm y)\in X_c(f)$, we have $f(\bm x)\leq g(\bm x,\bm y)$. When $m_2=0$ in Definition 1, $f(\bm x)$ is called a strong convertible convex(SCC) function on $S$.


When $S=R^n\times R^{m_1}\times R^{m_2}$, the term "on $S$" is omitted.

For each fixed $(\bm x,\bm y)\in R^n\times R^{m_1}$, let a set
\begin{eqnarray}
  X_c(\bm x,\bm y)=\{\bm z\mid (\bm x,\bm y,\bm z)\in X_c(f)\}. \label{d5}
\end{eqnarray}
For each fixed $\bm z\in S_3$, a set is defined by
\begin{eqnarray}
  Y_c(\bm z)=\{(\bm x,\bm y)\in S_1\times S_2\mid (\bm x,\bm y,\bm z)\in X_c(f)\}. \label{d6}
\end{eqnarray}

The following conclusion is clear.

\begin{prop}
Let $f$ be a SCN function on $S$. Then,
\begin{eqnarray}
\min\limits_{\bm x\in S_1} f(\bm x)=\min\limits_{(\bm x,\bm y,\bm z)\in X_c(f)} g(\bm x,\bm y,\bm z)=\min\limits_{(\bm x,\bm y)\in Y_c(\bm z)}\max\limits_{\bm z\in X_c(\bm x,\bm y)} g(\bm x,\bm y,\bm z)\label{d7}
\end{eqnarray}
and $-f$ is a SCCN function on $S$.
\end{prop}

Particularly, when $f$ is a SCC function, we have
\begin{eqnarray*}
\min\limits_{\bm x\in S_1} f(\bm x)=\min\limits_{(\bm x,\bm y)\in X_c(f)} g(\bm x,\bm y)=\min\limits_{(\bm x,\bm y)\in X_c(f)} g(\bm x,\bm y).
\end{eqnarray*}
Nonconvex optimization problem of the SCC function: $\min\limits_{\bm x\in S_1} f(\bm x)$ can be transformed into a convex optimization problem.

In fact, for a given $(\bar{\bm x},\bar{\bm y},\bar{\bm z})\in X(f)$, by \eqref{d4} we have
\begin{eqnarray*}
g(\bar{\bm x},\bar{\bm y},{\bm z})\leq f(\bar{\bm x})=g(\bar{\bm x},\bar{\bm y},\bar{\bm z})\leq g(\bar{\bm x},{\bm y},\bar{\bm z}), \forall (\bar{\bm x},\bar{\bm y},{\bm z}),(\bar{\bm x},{\bm y},\bar{\bm z})\in X_c(f).
\end{eqnarray*}
Because $X_c(f)$ is convex set, the right term of \eqref{d7} is rewritten as
\begin{eqnarray}
\min\limits_{\bm x\in S_1} f(\bm x)=\min\limits_{\bm x} \min\limits_{\bm y} \max\limits_{\bm z} g(\bm x,\bm y,\bm z)\ \ \mbox{s.t.}\ (\bm x,\bm y,\bm z)\in X_c(f) \label{d8}
\end{eqnarray}

When $(\bm x,\bm y,\bm z)\in X(f)$,  $(\bm x,\bm y,\bm z)$ is called a convertible nonconvex(SCN) point of $f$.

Next,  some  examples are given to show that the number of SCN form could be more than one.

\begin{ex}
Non-convex function $f(x_1,x_2)=2x_1x_2$ is a SCN function. One of its SCN forms is
\begin{eqnarray*}
g(x_1,x_2,z_1,z_2)&=&(x_1+x_2)^2-z_1-z_2:\\
g_1(x_1,x_2,z_1,z_2)&=& x_1^2-z_1,\\
g_2(x_1,x_2,z_1,z_2)&=&x_2^2-z_2.
\end{eqnarray*}
It is clear that $f(x_1,x_2)=g(x_1,x_2,z_1,z_2)$ for $(x_1,x_2,z_1,z_2)\in X(f)$ and
$$f(x_1,x_2)=\max\limits_{\bm z} g(x_1,x_2,z_1,z_2)\ \ \mbox{s.t.}\ (x_1,x_2,z_1,z_2)\in X_c(f).$$
A second SCN form of $f(x_1,x_2)=2x_1x_2$ is $[g(x_1,x_2,z_1)=(x_1+x_2)^2-z_1:g_1(x_1,x_2,z_1)=x_1^2+x_2^2-z_1]$.
Hence, it is understood that there are more than one SCN form.

If $f$ is a SCN function, $-f$ may also be a SCN function. In Example 2.1, $-f(x_1,x_2)=-2x_1x_2$ is a SCN function. One of its SCN forms is $[g(x_1,x_2,z_1)=(x_1-x_2)^2-z_1:g_1(x_1,x_2,z_1)=x_1^2+x_2^2-z_1]$.
\end{ex}%

\begin{ex}
(Example 2.1 in  Chen et al. (2014)) Nonsmooth function $f(x_1,x_2)=(x_1+x_2-1)^2+\lambda(|x_1|^\frac{1}{2}+|x_2|^\frac{1}{2})$ is converted to
\begin{eqnarray*}
g(x_1,x_2,y_1,y_2,y_3,y_4,z_1,z_2)&=&(x_1+x_2-1)^2+\lambda(y_1+y_3)\\
            && +y_1^4+x_1^2-2z_1+y_3^4+x_2^2-2z_2:\\
g_1(x_1,x_2,y_1,y_2,y_3,y_4,z_1,z_2)&=& y_1^4-z_1,\\
g_2(x_1,x_2,y_1,y_2,y_3,y_4,z_1,z_2)&=& x_1^2-z_1,\\
g_3(x_1,x_2,y_1,y_2,y_3,y_4,z_1,z_2)&=& y_2^2-y_1,\\
g_4(x_1,x_2,y_1,y_2,y_3,y_4,z_1,z_2)&=& y_3^4-z_2,\\
g_5(x_1,x_2,y_1,y_2,y_3,y_4,z_1,z_2)&=& x_2^2-z_2,\\
g_6(x_1,x_2,y_1,y_2,y_3,y_4,z_1,z_2)&=& y_4^2-y_3,
\end{eqnarray*}
where $\lambda>0$. So, $f(x)$ is a SCN function on $S=S_1\times S_2\times S_3$, where $S_1=R^2, S_2=\{\bm y\in R^4\mid y_1,y_3\in R^2_+, y_2,y_4\in R^1\}$ and $S_3=R_+^2$.
\end{ex}%

\begin{ex}
 Let the function $f(x_1,x_2)=(x_1+x_2-1)^2+\lambda\|(x_1,x_2)\|_0$ be nonconvex and discontinuous, where $\lambda>0$. Let $\bm x=(x_1,x_2),\bm y=(y_1,y_2),\bm z=(z_1,z_2)\in R^2$ and
\begin{eqnarray*}
g(\bm x,\bm y,\bm z)&=&(x_1+x_2-1)^2+\lambda(y_1+y_2)+(x_1+y_1-1)^2-z_1\\
&&+x_1^2+(y_1-1)^2-z_1+(x_2+y_2-1)^2-z_2+x_2^2+(y_2-1)^2-z_2:\\
g_1(\bm x,\bm y,\bm z)&=&(x_1+y_1-1)^2-z_1,\\
g_{2}(\bm x,\bm y,\bm z)&=&x_1^2+(y_1-1)^2-z_1,\\
g_{3}(\bm x,\bm y,\bm z)&=&y_1^2-y_1,\\
g_6(\bm x,\bm y,\bm z)&=&(x_2+y_2-1)^2-z_2,\\
g_{7}(\bm x,\bm y,\bm z)&=&x_2^2+(y_2-1)^2-z_2,\\
g_{9}(\bm x,\bm y,\bm z)&=&y_2^2-y_2.
\end{eqnarray*}
 So, $f(x)$ is a SCN function on $S=S_1\times S_2\times S_3$, where $S_1=R^2, S_2=\{\bm y\in R^2\mid 0\leq y_1,y_2\leq 1\}$ and $S_3=R_+^2$.
\end{ex}%
 The above examples show that many nonconvex and nonsmooth optimization problems can be solved by \eqref{d8} through its equivalent minmax problems.

\section{Operational Properties of SCN Function}

In this section, some operational properties of SCN function are proved. And it is always assumed that   SCN function is nonconvex or nondifferentiable.

\begin{prop}
If $f_1,f_2:R^n\to R$ are SCN functions  on convex $\bar{S}=S_1 \times \bar{S}_2\times \bar{S}_3 \subset R^n\times R^{m_1}\times R^{m_2}$ and $\tilde{S}={S}_1\times \tilde{S}_2\times\tilde{S}_3\subset R^n\times R^{m_3}\times R^{m_4}$ respectively, where $m_1,m_2, m_3,m_4\geq 0$, then
$\alpha_1 f_1+\alpha_2 f_2$ is a SCN function  on convex $S=S_1\times \bar{S}_2\times \tilde{S}_2 \times \bar{S}_3 \times\tilde{S}_3$ for any $\alpha_1,\alpha_2>0$.
\end{prop}%

{\it Proof.} Since $f_1$ and $f_2$ are SCN functions on $S$, their SCN forms  are given respectively by
   \begin{eqnarray}
f_1(\bm x)=[\bar{g}(\bm x,\bar{\bm y},\bar{\bm z}):\bar{g}_i(\bm x,\bar{\bm y},\bar{\bm z}), i=1,2,\cdots,s_1; \bar{h}_j(\bm x,\bar{\bm y},\bar{\bm z}), j=1,2,\cdots,r_1],\label{d9}\\
f_2(\bm x)=[\tilde{g}(\bm x,\tilde{\bm y},\tilde{\bm z}):\tilde{g}_i(\bm x,\tilde{\bm y},\tilde{\bm z}),  i=1,2,\cdots,s_2;\tilde{h}_j(\bm x,\tilde{\bm y},\tilde{\bm z}), j=1,2,\cdots,r_2], \label{d10}
\end{eqnarray}
where $\bar{g}_i(\bm x,\bar{\bm y},\bar{\bm z})$ is convex on $(\bm x,\bar{\bm y},\bar{\bm z})$,$i=1,2,\cdots,s_1$, $\bar{h}_j(\bm x,\bar{\bm y},\bar{\bm z})$ is linear on $(\bm x,\bar{\bm y},\bar{\bm z})$,$j=1,2,\cdots,r_1$,  $\tilde{g}_i(\bm x,\tilde{\bm y},\tilde{\bm z})$ is convex on $(\bm x,\tilde{\bm y},\tilde{\bm z})$,$i=1,2,\cdots,s_2$ and $\tilde{h}_j(\bm x,\tilde{\bm y},\tilde{\bm z})$ is linear on $(\bm x,\tilde{\bm y},\tilde{\bm z})$,$j=1,2,\cdots,r_2$. $\bar{g}(\bm x,\bar{\bm y},\bar{\bm z})$ and $\tilde{g}(\bm x,\tilde{\bm y},\tilde{\bm z})$ are convex on $(\bm x,\bar{\bm y})$ and $(\bm x,\tilde{\bm y})$ respectively. $\bar{g}(\bm x,\bar{\bm y},\bar{\bm z})$ and $\tilde{g}(\bm x,\tilde{\bm y},\tilde{\bm z})$ are concave on $\bar{\bm z}$ and $\tilde{\bm z}$ respectively. So, we have
\begin{eqnarray}
  X_c(f_1)=\{(\bm x,\bar{\bm y},\bar{\bm z})\in\bar{S} &\mid & \bar{g}_i(\bm x,\bar{\bm y},\bar{\bm z})\leq 0,\  i=1,2,\cdots,s_1,\nonumber\\
  && \bar{h}_j(\bm x,\bar{\bm y},\bar{\bm z})= 0,\  j=1,2,\cdots,r_1\}, \label{d12}
\end{eqnarray}
and
\begin{eqnarray}
  X_c(f_2)=\{(\bm x,\tilde{\bm y},\tilde{\bm z})\in \tilde{S}&\mid & \tilde{g}_i(\bm x,\tilde{\bm y},\tilde{\bm z})\leq 0,\  i=1,2,\cdots,s_2,\nonumber\\
   &&\tilde{h}_j(\bm x,\bar{\bm y},\bar{\bm z})= 0,\  j=1,2,\cdots,r_2\},\}.\label{d14}
\end{eqnarray}
Now, let $S=S_1\times \bar{S}_2\times \tilde{S}_2 \times \bar{S}_3 \times\tilde{S}_3$ and
\begin{eqnarray*}
X_c(\alpha_1 f_1+\alpha_2f_2)=\{(\bm x,\bar{\bm y},\tilde{\bm y},\bar{\bm z},\tilde{\bm z})\in S\mid (\bm x,\bar{\bm y},\bar{\bm z})\in X_c(f_1),(\bm x,\tilde{\bm y},\tilde{\bm z})\in X_c(f_2) \}.
\end{eqnarray*}
Since $f_1(\bm x)$ and $f_2(\bm x)$ are SCN functions , there are some $(\bm x,\bar{\bm y}^*,\bar{\bm z}^*)\in X(f_1)$ and  $(\bm x,\tilde{\bm y}^*,\tilde{\bm z}^*)\in X(f_2)$, such that
$$f_1(\bm x)=\bar{g}(\bm x,\bar{\bm y}^*,\bar{\bm z}^*) \mbox{ and } f_2(\bm x)=\tilde{g}(\bm x,\tilde{\bm y}^*,\tilde{\bm z}^*). $$

For each fixed $(\bm x,\bar{\bm y},\tilde{\bm y})$, let  sets
\begin{eqnarray*}
  X_c(\bm x,\bar{\bm y})=\{\bar{\bm z}\mid (\bm x,\bar{\bm y},\bar{\bm z})\in X_c(f_1)\},
  X_c(\bm x,\tilde{\bm y})=\{\tilde{\bm z}\mid (\bm x,\tilde{\bm y},\tilde{\bm z})\in X_c(f_2)\},
\end{eqnarray*}
and
\begin{eqnarray*}
  X_c(\bm x,\bar{\bm y},\tilde{\bm y})=\{(\bar{\bm z},\tilde{\bm z})\mid \bar{\bm z}\in X_c(\bm x,\bar{\bm y}),\tilde{\bm z}\in X_c(\bm x,\tilde{\bm y})\}.
\end{eqnarray*}
So, it is clear that we have
\begin{eqnarray*}
\alpha_1 f_1(\bm x)+\alpha_2 f_2(\bm x)&=&\alpha_1 \bar{g}(\bm x,\bar{\bm y}^*,\bar{\bm z}^*)+\alpha_2 \tilde{g}(\bm x,\tilde{\bm y}^*,\tilde{\bm z}^*)\\
&=&\alpha_1\min\limits_{\bar{\bm y}} \max\limits_{\bar{\bm z}} \bar{g}(\bm x,\bar{\bm y},\bar{\bm z})+\alpha_2\min\limits_{\tilde{\bm y}} \max\limits_{\tilde{\bm z}} \tilde{g}(\bm x,\tilde{\bm y},\tilde{\bm z})\\
&=&\min\limits_{(\bar{\bm y},\tilde{\bm y})} \max\limits_{(\bar{\bm z},\tilde{\bm z})}(\alpha_1 \bar{g}(\bm x,\bar{\bm y},\bar{\bm z})+\alpha_2 \tilde{g}(\bm x,\tilde{\bm y},\tilde{\bm z})).
\end{eqnarray*}
Hence, we have  $ (\bm x,\bar{\bm y}^*,\tilde{\bm y}^*,\bar{\bm z}^*,\tilde{\bm z}^*)\in X(\alpha_1 f_1+\alpha_2f_2)$. By the definition of SCN function, $\alpha_1 f_1+\alpha_2 f_2$ is a SCN function  on $S$.

\begin{prop}
Suppose that $f_1$ and $f_2$ are SCN functions  on $\bar{S}=S_1 \times \bar{S}_2\times \bar{S}_3 \subset R^n\times R^{m_1}\times R^{m_2}$ and $\tilde{S}={S}_1\times \tilde{S}_2\times\tilde{S}_3\subset R^n\times R^{m_3}\times R^{m_4}$ respectively, where $m_1, m_2, m_3, m_4\geq 0$, their SCN forms  are given respectively by \eqref{d9} and \eqref{d10}. If $\bar{g}(\bm x,\bar{\bm y},\bar{\bm z})\geq 0$ and $\tilde{g}(\bm x,\tilde{\bm y},\tilde{\bm z})\geq 0$ are convex on $(\bm x,\bar{\bm y},\bar{\bm z})\in \bar{S}$ and $(\bm x,\tilde{\bm y},\tilde{\bm z})\in \tilde{S}$ respectively, then
$f_1f_2$ is a SCN function  on $S=S_1\times \bar{S}_2\times \tilde{S}_2 \times \hat{S}_2 \times \bar{S}_3\times \tilde{S}_3\times \hat{S}_3$.
\end{prop}
{\it Proof.}  Since $f_1$ and $f_2$ are SCN functions on  $\bar{S}$ and $\tilde{S}$ respectively, their SCN forms are given respectively by \eqref{d9} and \eqref{d10}. Let $ \hat{\bm y}=(\hat{y}_1,\hat{y}_2)\in \hat{S}_2=R^2$,  $ \hat{\bm z}=\hat{z}_1\in \hat{S}_3=R^1$ and $S=S_1\times \bar{S}_2\times \tilde{S}_2 \times \hat{S}_2 \times \bar{S}_3\times \tilde{S}_3\times\hat{S}_3.$
A SCN form of $f_1(\bm x)f_2(\bm x)$ on $S$ be defined by
\begin{eqnarray*}
g(\bm x,\bar{\bm y},\tilde{\bm y},\hat{\bm y},\bar{\bm z},\tilde{\bm z},\hat{\bm z})&=&0.5(\hat{y}_1+\hat{y}_2)^2-0.5\hat{z}_1,\\
g_1(\bm x,\bar{\bm y},\tilde{\bm y},\hat{\bm y},\bar{\bm z},\tilde{\bm z},\hat{\bm z})&=&\hat{y}_1^2+\hat{y}_2^2-\hat{z}_1,\\
g_2(\bm x,\bar{\bm y},\tilde{\bm y},\hat{\bm y},\bar{\bm z},\tilde{\bm z},\hat{\bm z})&=&\bar{g}(\bm x,\bar{\bm y},\bar{\bm z})-\hat{y}_1,\\
g_3(\bm x,\bar{\bm y},\tilde{\bm y},\hat{\bm y},\bar{\bm z},\tilde{\bm z},\hat{\bm z})&=&\tilde{g}(\bm x,\tilde{\bm y},\tilde{\bm z})-\hat{y}_2,\\
g_{i+3}(\bm x,\bar{\bm y},\tilde{\bm y},\hat{\bm y},\bar{\bm z},\tilde{\bm z},\hat{\bm z})&=&\bar{g}_i(\bm x,\bar{\bm y},\bar{\bm z}),\  i=1,2,\cdots,s_1,\\
g_{j+r_1+3}(\bm x,\bar{\bm y},\tilde{\bm y},\hat{\bm y},\bar{\bm z},\tilde{\bm z},\hat{\bm z})&=&\tilde{g}_j(\bm x,\tilde{\bm y},\tilde{\bm z}),\  j=1,2,\cdots,s_2,\\
h_{j}(\bm x,\bar{\bm y},\tilde{\bm y},\hat{\bm y},\bar{\bm z},\tilde{\bm z},\hat{\bm z})&=&\bar{h}_j(\bm x,\bar{\bm y},\bar{\bm z}),\  j=1,2,\cdots,r_1,\\
h_{j+r_1}(\bm x,\bar{\bm y},\tilde{\bm y},\hat{\bm y},\bar{\bm z},\tilde{\bm z},\hat{\bm z})&=&\tilde{h}_j(\bm x,\tilde{\bm y},\tilde{\bm z}),\  j=1,2,\cdots,r_2,
\end{eqnarray*}
where $(\bm x,\bar{\bm y},\tilde{\bm y},\hat{\bm y},\bar{\bm z},\tilde{\bm z},\hat{\bm z})\in S$. Now, let
\begin{eqnarray*}
X_c(f_1f_2)=\{(\bm x,\bar{\bm y},\tilde{\bm y},\hat{\bm y},\bar{\bm z},\tilde{\bm z},\hat{\bm z})\in S&\mid&g_i(\bm x,\bar{\bm y},\tilde{\bm y},\hat{\bm y},\bar{\bm z},\tilde{\bm z},\hat{\bm z})\leq 0,i=1,2,3,\\
&& (\bm x,\bar{\bm y},\bar{\bm z})\in X_c(f_1),(\bm x,\tilde{\bm y},\tilde{\bm z})\in X_c(f_2) \}.
\end{eqnarray*}
It is clear that there is some $(\bm x,\bar{\bm y}^*,\bar{\bm z}^*)\in X(f_1),(\bm x,\tilde{\bm y}^*,\tilde{\bm z}^*)\in X(f_2)$ for fixed $\bm x\in S_1$ such that $f_1(\bm x)=\bar{g}(\bm x,\bar{\bm y}^*,\bar{\bm z}^*)$ and $f_2(\bm x)=\tilde{g}(\bm x,\tilde{\bm y}^*,\tilde{\bm z}^*)$. Then, let $\hat{y}_1^*=\bar{g}(\bm x,\bar{\bm y}^*,\bar{\bm z}^*)$,  $\hat{y}_2^*=\tilde{g}(\bm x,\tilde{\bm y}^*,\tilde{\bm z}^*)$ and $\hat{z}_1^*=\hat{y}_1^{*2}+\hat{y}_2^{*2}$. So, we have  $(\bm x,\bar{\bm y}^*,\tilde{\bm y}^*,\hat{\bm y}^*,\bar{\bm z}^*,\tilde{\bm z}^*,\hat{\bm z}^*)\in X(f_1f_2)$ and
\begin{eqnarray}
g(\bm x,\bar{\bm y}^*,\tilde{\bm y}^*,\hat{\bm y}^*,\bar{\bm z}^*,\tilde{\bm z}^*,\hat{\bm z}^*)=\hat{y}_1^*\hat{y}_2^*=\bar{g}(\bm x,\bar{\bm y}^*,\bar{\bm z}^*)\tilde{g}(\bm x,\tilde{\bm y}^*,\tilde{\bm z}^*)=f_1(\bm x)f_2(\bm x).\label{d15}
\end{eqnarray}
Next, let $(\bm x,\bar{\bm y},\tilde{\bm y},\hat{\bm y},\bar{\bm z},\tilde{\bm z},\hat{\bm z})\in X_c(f_1f_2)$ for fixed $\bm x$.
If $(\bm x,\bar{\bm y},\tilde{\bm y},\hat{\bm y},\bar{\bm z},\tilde{\bm z})$ is fixed, because $\hat{y}_1^2+\hat{y}_2^2\leq \hat{z}_1$ holds, we have
\begin{eqnarray}
\max\limits_{\hat{\bm z}} \{g(\bm x,\bar{\bm y},\tilde{\bm y},\hat{\bm y},\bar{\bm z},\tilde{\bm z},\hat{\bm z})&=&0.5(\hat{y}_1+\hat{y}_2)^2-0.5\hat{z}_1\}\nonumber\\
&=&\hat{y}_1\hat{y}_2.\label{d17}
\end{eqnarray}
So, for fixed $(\bm x,\bar{\bm y},\tilde{\bm y},\bar{\bm z},\tilde{\bm z})$,
because $0\leq \bar{g}(\bm x,\bar{\bm y},\bar{\bm z})\leq \hat{y}_1$ and $0\leq \tilde{g}(\bm x,\tilde{\bm y},\tilde{\bm z})\leq \hat{y}_2$, we have
\begin{eqnarray}
\min\limits_{\hat{\bm y}} \max\limits_{\hat{\bm z}} g(\bm x,\bar{\bm y},\tilde{\bm y},\hat{\bm y},\bar{\bm z},\tilde{\bm z},\hat{\bm z})&=&\min\limits_{\hat{\bm y}}\hat{y}_1\hat{y}_2\nonumber\\
&=& \bar{g}(\bm x,\bar{\bm y},\bar{\bm z}) \tilde{g}(\bm x,\tilde{\bm y},\tilde{\bm z}).\label{d18}
\end{eqnarray}
Hence, by \eqref{d17} and \eqref{d18} we obtian
\begin{eqnarray}
\min\limits_{(\bar{\bm y},\tilde{\bm y})}\max\limits_{(\bar{\bm z},\tilde{\bm z})}\min\limits_{\hat{\bm y}} \max\limits_{\hat{\bm z}} g(\bm x,\bar{\bm y},\tilde{\bm y},\hat{\bm y},\bar{\bm z},\tilde{\bm z},\hat{\bm z})&=&\min\limits_{(\bar{\bm y},\tilde{\bm y})}\max\limits_{(\bar{\bm z},\tilde{\bm z})}\min\limits_{\hat{\bm y}} \hat{y}_1\hat{y}_2\nonumber\\
&=&\min\limits_{(\bar{\bm y},\tilde{\bm y})} \max\limits_{(\bar{\bm z},\tilde{\bm z})} \bar{g}(\bm x,\bar{\bm y},\bar{\bm z}) \tilde{g}(\bm x,\tilde{\bm y},\tilde{\bm z})\nonumber\\
&=&\min\limits_{(\bar{\bm y},\tilde{\bm y})} \max\limits_{\bar{\bm z}} \bar{g}(\bm x,\bar{\bm y},\bar{\bm z}) \max\limits_{\tilde{\bm z}} \tilde{g}(\bm x,\tilde{\bm y},\tilde{\bm z})\nonumber\\
&=&\min\limits_{\bar{\bm y}} \max\limits_{\bar{\bm z}} \bar{g}(\bm x,\bar{\bm y},\bar{\bm z})\min\limits_{\tilde{\bm y}} \max\limits_{\tilde{\bm z}} \tilde{g}(\bm x,\tilde{\bm y},\tilde{\bm z})\nonumber\\
&=&f_1(\bm x)f_2(\bm x). \label{d19}
\end{eqnarray}
Now, for fixed $(\bar{\bm y},\tilde{\bm y})$ and $\hat{\bm z}$, we have
\begin{eqnarray}
\min\limits_{\hat{\bm y}}\max\limits_{(\bar{\bm z},\tilde{\bm z})}  g(\bm x,\bar{\bm y},\tilde{\bm y},\hat{\bm y},\bar{\bm z},\tilde{\bm z},\hat{\bm z})
=\max\limits_{(\bar{\bm z},\tilde{\bm z})}\min\limits_{\hat{\bm y}}g(\bm x,\bar{\bm y},\tilde{\bm y},\hat{\bm y},\bar{\bm z},\tilde{\bm z},\hat{\bm z}).\label{d20}
\end{eqnarray}
So, by \eqref{d15},\eqref{d19} and \eqref{d20} it is clear that
\begin{eqnarray*}
\min\limits_{(\bar{\bm y},\tilde{\bm y},\hat{\bm y})}\max\limits_{(\bar{\bm z},\tilde{\bm z},\hat{\bm z})} g(\bm x,\bar{\bm y},\tilde{\bm y},\hat{\bm y},\bar{\bm z},\tilde{\bm z},\hat{\bm z})&=&\min\limits_{(\bar{\bm y},\tilde{\bm y})}\min\limits_{\hat{\bm y}}\max\limits_{(\bar{\bm z},\tilde{\bm z})} \max\limits_{\hat{\bm z}} g(\bm x,\bar{\bm y},\tilde{\bm y},\hat{\bm y},\bar{\bm z},\tilde{\bm z},\hat{\bm z})\\
&=&\min\limits_{(\bar{\bm y},\tilde{\bm y})}\max\limits_{(\bar{\bm z},\tilde{\bm z})}\min\limits_{\hat{\bm y}} \max\limits_{\hat{\bm z}} g(\bm x,\bar{\bm y},\tilde{\bm y},\hat{\bm y},\bar{\bm z},\tilde{\bm z},\hat{\bm z})\\
&=&f_1(\bm x)f_2(\bm x).
\end{eqnarray*}
By the definition of SCN function, $f_1(\bm x)f_2(\bm x)$ is a SCN function  on $S$.

\begin{prop}
Suppose that $f$ is a SCCN function on $S=S_1\times S_2\times S_3$, a SCCN form of $-f$ is given by $[g:g_1,g_2,\cdots,g_s;h_1,h_2,\cdots,h_r]$. If $-g(\bm x,\bm y,\bm z)> 0$ is concave on $(\bm x,\bm y,\bm z)\in S$, then
$\frac{1}{f}$ is a SCN function on $\bar{S}=S_1\times \bar{S}_2\times \bar{S}_3$,
where $\bar{S}_2=\{(\bm y,\bar{y}_1,\bar{y}_2)\mid 0>\bar{y}_1>-\infty,0> \bar{y}_2>- \infty,\bm y\in S_2\}$ and $\bar{S_3}=\{\bar{z}_1\mid 0< \bar{z}_1<+\infty,\bm z\in S_3\}$.
\end{prop}

{\it Proof.}
Let $\bar{\bm y}=(\bar{y}_1,\bar{y}_2)$, $\bar{\bm z}=\bar{z}_1$ and
\begin{eqnarray*}
\bar{g}(\bm x,\bm y,\bar{\bm y}, \bm z,\bar{\bm z})& =&\bar{y}_1+(\bar{y}_1+\bar{y}_2)^2-\bar{z}_1-2+\bar{y}_1^2+\bar{y}_2^2-\bar{z}_1:\\
g_1(\bm x,\bm y,\bar{\bm y}, \bm z,\bar{\bm z})&=&(\bar{y}_1+\bar{y}_2)^2-\bar{z}_1-2,\\
 g_2(\bm x,\bm y,\bar{\bm y}, \bm z,\bar{\bm z})&=&\bar{y}_1^2+\bar{y}_2^2-\bar{z}_1,\\
 g_3(\bm x,\bm y,\bar{\bm y}, \bm z,\bar{\bm z})&=&g(\bm x,\bm y,\bm z)+\bar{y}_2,\\
  g_{3+i}(\bm x,\bm y,\bar{\bm y}, \bm z,\bar{\bm z})&=&g_i(\bm x,\bm y,\bm z), i=1,2,\cdots,s,\\
   h_{j}(\bm x,\bm y,\bar{\bm y}, \bm z,\bar{\bm z})&=&h_j(\bm x,\bm y,\bm z), j=1,2,\cdots,r.
  \end{eqnarray*}
Let $\bar{S}=S_1\times \bar{S}_2\times \bar{S}_3$ and
\begin{eqnarray*}
X_c(\frac{1}{f})=\{(\bm x,\bm y,\bar{\bm y}, \bm z,\bar{\bm z})\in \bar{S}\mid g_i(\bm x,\bm y,\bar{\bm y}, \bm z,\bar{\bm z})\leq 0,i=1,2,3,
 (\bm x,\bm y,\bm z)\in X_c(f) \}.
\end{eqnarray*}
It is clear that $\bar{g}(\bm x,\bm y^*,\bar{\bm y}^*, \bm z^*,\bar{\bm z}^*) =\bar{y}_1^*=\frac{1}{-g(\bm x,\bm y^*,\bm z^*)}=\frac{1}{f(\bm x)}$ because there is some  $(\bm x,\bm y^*,\bar{\bm y}^*, \bm z^*,\bar{\bm z}^*)\in X(\frac{1}{f})$.
For $(\bm x,\bm y,\bar{\bm y}, \bm z,\bar{\bm z})\in X_c(\frac{1}{f})$, we have
\begin{eqnarray}
\min\limits_{\bar{\bm y}}\max\limits_{\bar{\bm z}}\bar{g}(\bm x,\bm y,\bar{\bm y}, \bm z,\bar{\bm z})& =&\min\limits_{\bar{\bm y}}\max\limits_{\bar{\bm z}}\bar{y}_1+(\bar{y}_1+\bar{y}_2)^2-\bar{z}_1-2+\bar{y}_1^2+\bar{y}_2^2-\bar{z}_1 \nonumber\\
&=&\min\limits_{\bar{\bm y}} \bar{y}_1\nonumber\\
&=&\min\limits_{\bar{\bm y}}\frac{1}{\bar{y}_2}\nonumber\\
&=& \frac{1}{-g(\bm x,\bm y,\bm z)}, \label{d201}
\end{eqnarray}
where $\bar{y}_1\bar{y}_2=1$ and $\bar{y}_2\leq  -g(\bm x,\bm y,\bm z) $ for $(\bm x,\bm y,\bm z)\in X_c(f)$. For $(\bm x,\bm y,\bm z)\in X_c(f)$, we have
$$\min\limits_{\bm y}  \max\limits_{\bm z}  g(\bm x,\bm y,\bm z)=-f(\bm x). $$
Hence, by \eqref{d201}, we have
\begin{eqnarray*}
\min\limits_{(\bm y,\bar{\bm y})}\max\limits_{(\bm z,\bar{\bm z})}\bar{g}(\bm x,\bm y,\bar{\bm y}, \bm z,\bar{\bm z}))
&=&\min\limits_{\bm y}\max\limits_{\bm z}\min\limits_{\bar{\bm y}}\max\limits_{\bar{\bm z}}\bar{g}(\bm x,\bm y,\bar{\bm y}, \bm z,\bar{\bm z})\\
&=& \min\limits_{\bm y}\max\limits_{\bm z}\frac{1}{-g(\bm x,\bm y,\bm z)},\\
&=& \frac{1}{\max\limits_{\bm y}\min\limits_{\bm z}(-g(\bm x,\bm y,\bm z))},\\
&=& -\frac{1}{\min\limits_{\bm y}\max\limits_{\bm z}g(\bm x,\bm y,\bm z)},\\
&=& \frac{1}{f(\bm x)}.
\end{eqnarray*}
So, $\frac{1}{f(\bm x)}$ is a SCN function.

By  Proposition 3.3, we have the following conclusion.

\begin{prop}
Suppose that $f$ is a SCN function on $S=S_1\times S_2\times S_3$, a SCN form of $f$ is given by $[g:g_1,g_2,\cdots,g_s;h_1,h_2,\cdots,h_r]$. If $-g(\bm x,\bm y,\bm z)< 0$ is concave on $(\bm x,\bm y,\bm z)\in S$, then
$\frac{1}{f}$ is a SCCN function on $\bar{S}=S_1\times \bar{S}_2\times \bar{S}_3$,
where $\bar{S}_2=\{(\bm y,\bar{y}_1,\bar{y}_2)\mid 0>\bar{y}_1>-\infty,0> \bar{y}_2>- \infty,\bm y\in S_2\}$ and $\bar{S_3}=\{\bar{z}_1\mid 0< \bar{z}_1<+\infty,\bm z\in S_3\}$.
\end{prop}
\begin{prop} 
  If $f:R^n\to R$ is a SCN function  on $S$ and $\phi:R\to R$ is a monotone increasing convex function, then $\phi(f(\bm x))$ is a SCN function  on $S$.
\end{prop}

{\it Proof.} Since $f(\bm x)$ is a SCN function  on $S$, a SCN form of $\phi(f(\bm x))$ is given by
$f=[\phi(g):g_1,g_2,$ $\cdots,g_s;h_1,h_2,\cdots,h_r]$.
Now, let $X_c(\phi(f))=X_c(f)$.
So, by Definition 1, there is some $(\bm x,\bm y,\bm z)\in X(\phi(f))$ such that
$$\phi(f(\bm x))=\phi(g(\bm x,\bm y,\bm z)), $$
 Hence, we have
\begin{eqnarray*}
 \phi(f(\bm x))=\min\limits_{\bm y} \max\limits_{\bm z} \phi(g(\bm x,\bm y,\bm z))\ \ \mbox{s.t.}\ (\bm x,\bm y,\bm z)\in X_c(f).
\end{eqnarray*}
By the definition of SCN function, $\phi(f)$ is a  SCN function  on $S$.

\begin{prop}
 Suppose that $f$ is a SCN function on $S=S_1\times S_2\times S_3$, a SCN form of $f$ is given by $[g:g_1,g_2,\cdots,g_s;h_1,h_2,\cdots,h_r]$. If $g(\bm x,\bm y,\bm z)\geq 0$ is convex on $(\bm x,\bm y,\bm z)\in S$, then
$f(\bm x)^a$ is a SCN function on $\bar{S}=S_1\times \bar{S}_2\times \bar{S}_3$,
where $a>0$, $\bar{S}_2=\{(\bm y,\bar{y}_1,\bar{y}_2)\mid 0\leq \bar{y}_1<+\infty,0\leq \bar{y}_2\leq +\infty,\bm y\in S_2\}$ and $\bar{S_3}=\{\bar{z}_1\mid 0\leq \bar{z}_1<+\infty,\bm z\in S_3\}$.
\end{prop}

{\it Proof.}  When $0<a<1$, let $\bar{\bm y}=(\bar{y}_1,\bar{y}_2)$, $\bar{\bm z}=\bar{z}_1$ and
\begin{eqnarray*}
\bar{g}(\bm x,\bm y,\bar{\bm y}, \bm z,\bar{\bm z})& =&\bar{y}_1+(\bar{y}_1)^\frac{2}{a}-\bar{z}_1+\bar{y}_2^2-\bar{z}_1:\\
g_1(\bm x,\bm y,\bar{\bm y}, \bm z,\bar{\bm z})&=&(\bar{y}_1)^\frac{2}{a}-\bar{z}_1,\\
g_2(\bm x,\bm y,\bar{\bm y}, \bm z,\bar{\bm z})&=&\bar{y}_2^2-\bar{z}_1,\\
g_3(\bm x,\bm y,\bar{\bm y}, \bm z,\bar{\bm z})&=&g(\bm x,\bm y,\bm z)-\bar{y}_2,\\
g_{3+i}(\bm x,\bm y,\bar{\bm y}, \bm z,\bar{\bm z})&=&g_i(\bm x,\bm y,\bm z), i=1,2,\cdots,s,\\
  h_{j}(\bm x,\bm y,\bar{\bm y}, \bm z,\bar{\bm z})&=&h_j(\bm x,\bm y,\bm z), j=1,2,\cdots,r.
  \end{eqnarray*}
Let $\bar{S}=S_1\times \bar{S}_2\times \bar{S}_3$ and
\begin{eqnarray*}
X_c(f^a)=\{(\bm x,\bm y,\bar{\bm y}, \bm z,\bar{\bm z})\in \bar(S)\mid g_i(\bm x,\bm y,\bar{\bm y}, \bm z,\bar{\bm z})\leq 0,i=1,2,3,
 (\bm x,\bm y,\bm z)\in X_c(f) \}.
\end{eqnarray*}
It is clear that $\bar{g}(\bm x,\bm y^*,\bar{\bm y}^*, \bm z^*,\bar{\bm z}^*) =\bar{y}_1^*=\bar{y}_2^{*a}={f(\bm x)}^a$ because there is some  $(\bm x,\bm y^*,\bar{\bm y}^*, \bm z^*,\bar{\bm z}^*)\in X({f}^a)$.
For $(\bm x,\bm y,\bar{\bm y}, \bm z,\bar{\bm z})\in X_c({f}^a)$, we have
\begin{eqnarray*}
\min\limits_{\bar{\bm y}}\max\limits_{\bar{\bm z}}\bar{g}(\bm x,\bm y,\bar{\bm y}, \bm z,\bar{\bm z})& =&\min\limits_{\bar{\bm y}}\max\limits_{\bar{\bm z}}\bar{y}_1+(\bar{y}_1)^\frac{2}{a}-\bar{z}_1+\bar{y}_2^2-\bar{z}_1 \nonumber\\
&=&\min\limits_{\bar{\bm y}} \bar{y}_1\nonumber\\
&=& g(\bm x,\bm y,\bm z)^a,
\end{eqnarray*}
where $\bar{y}_1\bar{y}_2=1$ and $g(\bm x,\bm y,\bm z)\leq \bar{y}_2  $ for $(\bm x,\bm y,\bm z)\in X_c(f)$. For $(\bm x,\bm y,\bm z)\in X_c(f)$, we have
$$\min\limits_{\bm y}  \max\limits_{\bm z}  g(\bm x,\bm y,\bm z)= f(\bm x). $$
Hence, we have
\begin{eqnarray*}
\min\limits_{(\bm y,\bar{\bm y})}\max\limits_{(\bm z,\bar{\bm z})}\bar{g}(\bm x,\bm y,\bar{\bm y}, \bm z,\bar{\bm z}))
&=&\min\limits_{\bm y}\max\limits_{\bm z}\min\limits_{\bar{\bm y}}\max\limits_{\bar{\bm z}}\bar{g}(\bm x,\bm y,\bar{\bm y}, \bm z,\bar{\bm z})\\
&=& \min\limits_{\bm y}\max\limits_{\bm z}g(\bm x,\bm y,\bm z)^a,\\
&=& {f(\bm x)}^a.
\end{eqnarray*}
So, ${f(\bm x)}^a$ is a SCN function on $S$.

 When $1<a<2$, ${f(\bm x)}^a=f(\bm x)f(\bm x)^{a-1}$, then ${f(\bm x)}$ and  $f(\bm x)^{a-1}$ are SCN functions on $S$ and $\bar{S}$ respectively. By Proposition 3.2, ${f(\bm x)}^a$ is a SCN function on $\bar{S}$.

 When $a\geq 2$, ${f(\bm x)}^a$ is a SCN function on $S$  by  Proposition 3.2.

\begin{cor}
 Suppose that $f_1(\bm x)$ and $f_2(\bm x)$ are SCN functions  on $\bar{S}=S_1 \times \bar{S}_2\times \bar{S}_3 \subset R^n\times R^{m_1}\times R^{m_2}$ and $\tilde{S}={S}_1\times \tilde{S}_2\times\tilde{S}_3\subset R^n\times R^{m_3}\times R^{m_4}$, where $m_1, m_2, m_3, m_4\geq 0$, their SCN forms  are given respectively by \eqref{d9} and \eqref{d10}. Let $a,b>0$. If $\bar{g}(\bm x,\bar{\bm y},\bar{\bm z})\geq 0$ and $\tilde{g}(\bm x,\tilde{\bm y},\tilde{\bm z})\geq 0$ are convex on $(\bm x,\bar{\bm y},\bar{\bm z})\in \bar{S}$ and $(\bm x,\tilde{\bm y},\tilde{\bm z})\in \tilde{S}$ respectively, then
$f_1(\bm x)^af_2(\bm x)^b$ is a SCN function  on $S=S_1\times \bar{S}_2\times \tilde{S}_2 \times \hat{S}_2 \times \bar{S}_3\times \tilde{S}_3\times \hat{S}_3$.
\end{cor}

By Proposition 3.1-3.6, some polynomial functions are SCN functions. For Example, multi-convex function $f(\bm x)=x_1x_2\cdots x_n$ is a SCN function. Therefore,  SCN  function covers a wide range of non-convex functions.
The following function are SCN functions.


\begin{ex}
(1) $f(\bm x)=\frac{a_1^\top\bm x+b_1}{a_2^\top\bm x+b_2}$ is a SCN function,
 where $a_1,a_2\in R^n$ are given and $\bm x\in R^n$ is variable. 
(2) $f(\bm x)$ is a SCN function, where $f(\bm x) $ be a concave function on $S_1$,

(3)  $f(\bm x)=\frac{1}{b(\bm x)}$ a SCN function, where $b(\bm x)$ is concave  on convex $S_1\subset R^n$ and $b(\bm x)>0$.

(4) $f(x)=\frac{-1}{1+exp(x)}$ on $x\in R^1$ is a SCN function.

(5) $f(\bm x)=\frac{1}{b(\bm x)}$ is a SCN function, where $b(\bm x)$ is a convex function on convex $S_1\subset R^n$ and $b(\bm x)<0$.

(6) $f(\bm x)=f_1(\bm x)f_2(\bm x)$ is a SCN functio, where $f_1(\bm x),f_2(\bm x)$ are convex  on convex set $S_1\subset R^n$ and $f_1(\bm x)\geq 0$ and $f_2(\bm x)\geq 0$ are true on $S_1$.

(7)If $f_1(\bm x)$ and $f_2(\bm x)$ are bounded and convex  on convex set $S_1\subset R^n$ and $f_1(\bm x)\geq M,f_2(\bm x)\geq N$ for given $N,M<0$, then $f(\bm x)=f_1(\bm x)f_2(\bm x)=(f_1(\bm x)-M)(f_2(\bm x)-N)+Nf_1(\bm x)+Mf_2(\bm x)-NM$ is a SCN function on $S$,
 where  $S_2=\{(y_1,y_2)\mid M\leq y_1<+\infty, N\leq  y_2\leq +\infty\}$, $S_3=\{(z_1,z_2,z_3)\mid 0\leq z_1<+\infty,M\leq z_2<+\infty,N\leq z_3<+\infty\}$ and $S=S_1\times S_2 \times S_3$.

\end{ex}%




The SCN forms of some functions: trigonometric function, DC function, entropy function, 0-norm function, sigmoid function are given, symbolic function and so on as follows.

\begin{ex}
(1) For functions like $\sin x$, different $S_1\in R^1$ brings different SCN forms as follows.

$\bullet$ On $x\in S_1=[0,\pi]$, a SCN form of function $f(x)=\sin x$ is $[z_1:-\sin x+z_1]$, where $S=S_1\times S_3$ and $S_3=\{z_1\mid 0\leq z_1\leq 1\}.$

$\bullet$ On $x\in S_1=[-\pi,\pi]$, a CN form of function $f(x)=\sin x$ is $[y_1-\sin (\frac{\pi+x}{2})-z_1-z_1^2-z_2^2+1:(z_1+z_2)^2-1-y_1,-\sin (\frac{\pi+x}{2})-z_1,z_1^2+z_2^2-1]$, where $S=S_1\times S_2\times S_3$, $S_2=\{y_1\mid -1\leq y_1\leq 1 \}$ and $S_3=\{(z_1,z_2)\mid -1\leq z_1\leq 1,-1\leq z_2\leq 1\}$.

$\bullet$ On $x\in S=[0,2\pi]$, a CN form of function $f(x)=\sin x$ is $[y_1-\sin \frac{x}{2}+z_1-z_1^2-z_2^2+1:(z_1+z_2)^2-1-y_1,-\sin \frac{x}{2}+z_1,z_1^2+z_2^2-1]$, where $S=S_1\times S_2\times S_3$, $S_2=\{y_1\mid -1\leq y_1\leq 1 \}$ and $S_3=\{(z_1,z_2)\mid -1\leq z_1\leq 1,-1\leq z_2\leq 1\}$.

$\bullet$ And on $x\in S=[0,2\pi]$, a CN form of function $f(x)=\cos x$ is $[z_3-(z_1+z_2)^2+1+z_4-\sin \frac{x}{2}+z_1-z_1^2-z_2^2+1-z_3^2-z_4^2+1:(z_1+z_2)^2-1-z_4,-\sin \frac{x}{2}+z_1,z_1^2+z_2^2-1,z_3^2+z_4^2-1]$, where $S=S_1\times S_3$ and $S_3=\{(z_1,z_2,z_3,z_4)\mid -1\leq z_i\leq 1,i=1,2,3,4\}$.

$\bullet$ Let $f(x)=\sin x$ on $x\in S=[b,a]$ with $b-a>0$. For any $x\in [a,b]$, let $K=[\frac{b-a}{2\pi}]$. There are $k\in \{0,1,2\cdots,K\}$ and $x'\in [0,2\pi]$ such that $x-a=2k\pi+x'$ and
 $$\sin x=\sin (2k\pi+x'+a)=\sin x'\cos a+\cos x' \sin a. $$
 So, when $\sin a=0$ and $x'\in [0,2\pi]$, we have a CN form:
 $$\sin x=[y_1-(-1)^{[\frac{x-a}{2\pi}]}\sin \frac{x-a}{2}+z_1-z_1^2-z_2^2+1:(z_1+z_2)^2-1-y_1,-(-1)^{[\frac{x-a}{2\pi}]}\sin \frac{x-a}{2}+z_1,z_1^2+z_2^2-1],$$
 where $S=S_1\times S_2\times S_3$, $S_2=\{y_1\mid -1\leq y_1\leq 1 \}$ and $S_3=\{(z_1,z_2)\mid -1\leq z_1\leq 1,-1\leq z_2\leq 1\}$.
When $\sin a\not=0$ and $x'\in [0,2\pi]$, we have a CN form:
\begin{eqnarray*}
 \sin x=[-(z_1+z_2)^2+1+z_4-(-1)^{[\frac{x-a}{2\pi}]}\sin \frac{x-a}{2}+z_1-z_1^2-z_2^2+1-z_3^2-z_4^2+1+z_4\cos a+z_3 \sin a:\\
 (z_1+z_2)^2-1-z_4,-(-1)^{[\frac{x-a}{2\pi}]}\sin \frac{x-a}{2}+z_1,z_1^2+z_2^2-1,z_3^2+z_4^2-1],
\end{eqnarray*}
where $S=S_1\times S_3$ and $S_3=\{(z_1,z_2,z_3,z_4)\mid -1\leq z_i\leq 1,i=1,2,3,4\}$.

(2) DC function is a very important class of nonconvex functions in THi et al(2018). Since the DC function $f(\bm x)=d(\bm x)-c(\bm x)$ can be converted to
$$f(\bm x)=[g(\bm x,z)=d(\bm x)-z: g_1(\bm x,z)=c(\bm x)-z],$$
where $d(\bm x)$ and $c(\bm x)$ are convex on $\bm x$, it is a SCN function.
By Hartman (1959), if $f:R^n\to R^1$ is a second-order continuously differentiable function on $R^n$, $f(x)$ is a DC function. Hence, $f(x)$ is a SCN function. \underline{ Hence, all second-order continuously differentiable functions on $R^n$ are SCN} \underline{functions.}

(3) A SCN form of entropy function $f(\bm x)=-\sum\limits_{i=1}^n x_i\ln x_i$ on $S_1=\{\bm x\in R^n\mid 0<x_i\leq 1\}$
 is defined by
 \begin{eqnarray*}
[\sum_{i=1}^m (0.5(x_i+y_i)^2-0.5z_i):
-\ln x_i-y_i, x_i^2+y_i^2-z_i; i=1,2,\cdots,n],
\end{eqnarray*}
 where $S_2=\{\bm y\in R^n\mid 0\leq y_i<+\infty,i=1,2,\cdots,n\}$ and $S_3=\{\bm z\in R^n\mid 0\leq z_i<+\infty,i=1,2,\cdots,n\}$.

(4) A SCN form of 0-norm function  $f(\bm x)=\|\bm x\|_0$ is a case of Example 3.5 when $A\bm x=b$.

(5)  A SCN form of sigmoid function $f(x)=\frac{2}{1+exp(-x)}-1$ is defined by
 \begin{eqnarray*}
[2y_1-1+(y_1+y_2)^2-z_1-1+y_1^2+y_2^2-z_1:(y_1+y_2)^2-z_1-1,y_1^2+y_2^2-z_1,1+exp(-x)-y_2],
\end{eqnarray*}
where $S_1=R^1,S_2=\{(y_1,y_2)^\top\mid 1\leq y_1\leq +\infty,1\leq y_2\leq +\infty\},S_3=\{z_1\mid 0\leq z_1\leq +\infty\}$.

(6) For power functions like $x^a$, different $S_1\in R^1$ and $a\in R^1$ brings different SCN forms as follows.

$\bullet$ On $x\in S_1=R^1$ and $0<a< 1$, a SCN form of function $f(x)=x^a$ is $[y_1+(y_1^{\frac{1}{a}})^2-z_1+x^2-z_1:(y_1^{\frac{1}{a}})^2-z_1,x^2-z_1]$, where $S=S_1\times S_2\times S_3$, $S_2=R^1$ and $S_3=R_+^1$.

$\bullet$ On $x\in S_1=R^1_+$ and $0<a< 1$, a SCN form of function $f(x)=x^a$ is $[y_1+(y_1^{\frac{1}{a}})^2-z_1+x^2-z_1:(y_1^{\frac{1}{a}})^2-z_1,x^2-z_1]$, where $S=S_1\times S_2\times S_3$, $S_2=R^1_+$ and $S_3=R^1_+$.

$\bullet$ On $x\in S_1=R^1$ and $1<a+1< 2$, a SCN form of function $f(x)=x^{a+1}$ is $[0.5(y_1+x)^2-0.5z_1-0.5z_2:(y_1^{\frac{1}{a}})^2-z_1,x^2-z_1,y_1^2-z_2]$, where $S=S_1\times S_2\times S_3$, $S_2=R^1$ and $S_3=R^2_+$.

$\bullet$ On $x\in S_1=R^1_+$ and $1<a+1< 2$, a SCN form of function $f(x)=x^{a+1}$ is $[0.5(y_1+x)^2-0.5z_1-0.5z_2:(y_1^{\frac{1}{a}})^2-z_1,x^2-z_1,y_1^2-z_2]$, where $S=S_1\times S_2\times S_3$, $S_2=R^1_+$ and $S_3=R^2_+$.

$\bullet$ On $x\in S_1=R^1$, $0<a<1$ and $n\geq 1$,  a SCN form of function $f(x)=x^{a+2n}$ is $[0.5(y_1+y_2)^2-0.5z_1+(y_1^{\frac{1}{a}})^2-z_2+x^2-z_2:(y_1^{\frac{1}{a}})^2-z_2$, $x^2-z_2,x^{2n}-y_2,y_1^2+y_2^2-z_1]$, where $S=S_1\times S_2\times S_3$, $S_2=R^1\times R^1_+$ and $S_3=R^2_+$.

$\bullet$ On $x\in S_1=R^1_+$, $0<a<1$ and $n\geq 1$,  a SCN form of function $f(x)=x^{a+2n}$ is $[0.5(y_1+y_2)^2-0.5z_1+(y_1^{\frac{1}{a}})^2-z_2+x^2-z_2:(y_1^{\frac{1}{a}})^2-z_2$, $x^2-z_2,x^{2n}-y_2,y_1^2+y_2^2-z_1]$, where $S=S_1\times S_2\times S_3$, $S_2=R^2_+$ and $S_3=R^2_+$.

(7) For the symbolic function function $f(x)=sgn(x)=$
\begin{math} \left\{
\begin{array}{ll}
1,& if \kg x> 0,\\
0,& if \kg x= 0,\\
-1,& if \kg x< 0,\\
\end{array}\right.\end{math}
a SCN form of $f(x)=sgn(x)$ is defined by
 \begin{eqnarray*}
g(x,\bm y,\bm z)=y_1+y_1^2-1+y_4+\sum\limits_{i=1}^9 g_i(x,\bm y,\bm z):&&\\
g_1(x,\bm y,\bm z)=(y_2+y_1-1)^2-z_1,&&
g_2(x,\bm y,\bm z)= y_2^2+(y_1-1)^2-z_1, \\
g_3(x,\bm y,\bm z)=(y_3+y_1)^2-z_2 ,&&
g_4(x,\bm y,\bm z)= y_3^2+(y_1)^2-z_2, \\
g_5(x,\bm y,\bm z)= y_2^2-x-z_3, &&
g_6(x,\bm y,\bm z)= y_3^2-z_3, \\
g_7(x,\bm y,\bm z)=(y_4+y_1)^2-z_4 ,&&
g_8(x,\bm y,\bm z)= y_4^2+y_1^2-z_4, \\
g_9(x,\bm y,\bm z)= y_1^2-z_5,&&
h_{1}(x,\bm y,\bm z)= z_5-1+y_4,
\end{eqnarray*}
where $\bm y\in R^4,\bm z\in R^5$, $S_1=R^1,S_2=R^4,S_3=R^5_+$.

Another SCN form of symbolic function function $f(x)=sgn(x)$ is defined by
 \begin{eqnarray*}
g(x,\bm y,\bm z)=y_1+y_1^2-1+y_4+\sum\limits_{i=1}^7 g_i(x,\bm y,\bm z):&&\\
g_1(x,\bm y,\bm z)=y_1^2-z_1,&&
g_2(x,\bm y,\bm z)=y_2^2-z_2, \\
g_3(x,\bm y,\bm z)=y_3^2-z_3,&&
g_4(x,\bm y,\bm z)=y_4^2-z_4, \\
g_5(x,\bm y,\bm z)=(y_2+y_1-1)^2-z_5,&&
g_6(x,\bm y,\bm z)=(y_3+y_1+1)^2-z_6, \\
g_7(x,\bm y,\bm z)=(y_4+y_1)^2-z_7 ,&&
h_1(x,\bm y,\bm z)= z_2+z_1-2y_1+1-z_5, \\
h_2(x,\bm y,\bm z)=z_3+z_1+2y_1+1-z_6,&&
h_{3}(x,\bm y,\bm z)= z_2-x-z_3, \\
h_{4}(x,\bm y,\bm z)= z_4+z_1-z_7,&&
h_{5}(x,\bm y,\bm z)= z_1-1+y_4,
\end{eqnarray*}
where $\bm y\in R^4,\bm z\in R^7$, $S_1=R^1,S_2=R^4,S_3=R^7_+$.

(8) In machine learning, there is a symbolic function function $f(x)=sgn(x)=$
\begin{math} \left\{
\begin{array}{ll}
1,& if \kg x\geq 0,\\
0,& if \kg x< 0,\\
\end{array}\right.\end{math}
then a SCN form of $f(x)=sgn(x)$ is defined by
 \begin{eqnarray*}
g(x,\bm y,\bm z)=y_1+\sum\limits_{i=1}^6 g_i(x,\bm y,\bm z):&&\\
g_1(x,\bm y,\bm z)=(y_2+y_1-1)^2-z_1 , &&
g_2(x,\bm y,\bm z)= y_2^2+(y_1-1)^2-z_1, \\
g_3(x,\bm y,\bm z)=(y_3+y_1)^2-z_2 , &&
g_4(x,\bm y,\bm z)= y_3^2+y_1^2-z_2, \\
g_5(x,\bm y,\bm z)= y_2^2-x-z_3, &&
g_6(x,\bm y,\bm z)= y_3^2-z_3,
\end{eqnarray*}
where $\bm y\in R^3,\bm z\in R^3$, $S_1=R^1,S_2=[0,1]\times R^2,S_3=R^3_+$.

Another SCN form of symbolic function function $f(x)=sgn(x)$ is defined by
 \begin{eqnarray*}
g(x,\bm y,\bm z)=y_1+\sum\limits_{i=1}^5 g_i(x,\bm y,\bm z):&&\\
g_1(x,\bm y,\bm z)=y_1^2-z_1, &&
g_2(x,\bm y,\bm z)=y_2^2-z_2, \\
g_3(x,\bm y,\bm z)=y_3^2-z_3, &&
g_4(x,\bm y,\bm z)=(y_2+y_1-1)^2-z_4, \\
g_5(x,\bm y,\bm z)=(y_3+y_1)^2-z_5 , &&
h_1(x,\bm y,\bm z)= z_2+z_1-2y_1+1-z_4, \\
h_2(x,\bm y,\bm z)= z_3+z_1-z_5,&&
h_3(x,\bm y,\bm z)= z_2-x-z_3,
\end{eqnarray*}
where $\bm y\in R^3,\bm z\in R^5$, $S_1=R^1,S_2=[0,1]\times R^2,S_3=R^5_+$.

(9) A SCN form of $f(x)=\max\{x,0\}$ is defined by
 \begin{eqnarray*}
g(x,\bm y,\bm z)=y_1^2+\sum\limits_{i=1}^4 g_i(x,\bm y,\bm z):&&\\
g_1(x,\bm y,\bm z)=(y_2+y_1)^2-z_1 , &&
g_2(x,\bm y,\bm z)= y_2^2+y_1^2-z_1, \\
g_3(x,\bm y,\bm z)= y_1^2-x-z_2, &&
g_4(x,\bm y,\bm z)= y_2^2-z_2,
\end{eqnarray*}
where $\bm y\in R^2,\bm z\in R^2$, $S_1=R^1,S_2=R^2,S_3=R^2_+$.
Another SCN form of $f(x)=\max\{x,0\}$ is defined by
 \begin{eqnarray*}
g(x,\bm y,\bm z)=y_1^2+\sum\limits_{i=1}^3 g_i(x,\bm y,\bm z):&&\\
g_1(x,\bm y,\bm z)=y_1^2-z_1, &&
g_2(x,\bm y,\bm z)=y_2^2-z_2, \\
g_3(x,\bm y,\bm z)=(y_2+y_1)^2-z_3, &&
h_1(x,\bm y,\bm z)= z_1+z_2-z_3, \\
h_2(x,\bm y,\bm z)= z_1-x-z_2,&&
\end{eqnarray*}
where $\bm y\in R^2,\bm z\in R^3$, $S_1=R^1,S_2=R^2,S_3=R^3_+$.

(10) Let $b(\bm x)$ be convex on $R^n$. So, its SCN form of $f(x)=\max\{b(\bm x),0\}$ is defined by
 \begin{eqnarray*}
g(\bm x,\bm y,\bm z)=y_1^2+\sum\limits_{i=1}^3 g_i(x,\bm y,\bm z):&&\\
g_1(\bm x,\bm y,\bm z)=y_1^2-z_1, &&
g_2(\bm x,\bm y,\bm z)=y_2^2-z_2, \\
g_3(\bm x,\bm y,\bm z)=(y_2+y_1)^2-z_3, &&
g_4(\bm x,\bm y,\bm z)=b(\bm x)-z_4, \\
h_1(\bm x,\bm y,\bm z)= z_1+z_2-z_3, &&
h_2(\bm x,\bm y,\bm z)= z_1-z_4-z_2,
\end{eqnarray*}
where $\bm y\in R^2,\bm z\in R^4$, $S_1=R^n,S_2=R^2,S_3=R^4_+$.
\end{ex}

It is easy to obtain SCN forms of some special functions, such as $\max\{b(\bm x),0\}^a(0<a<1)$ and $\max\{b_1(\bm x),b_2(\bm x)\}$.

Next, let us see more complex SCN forms of some special structure functions.

\begin{ex}
Let a polynomial function optimization be defined by
\begin{eqnarray}
\min\limits_{\bm x\in R^n} \Phi_1(\bm x)=\sum_{i=1}^m a_i\prod\limits_{j=1}^n x_j^{\alpha_{ij}},\label{d21}
\end{eqnarray}
where $a_i\not=0,\alpha_{ij}\in R^1,  i=1,2,\cdots,m, j=1,2,\cdots,n, S_1=\{\bm x\in R^n| x_i>0,i=1,2,\cdots,n\}.$ Let $\bm x=(x_1,x_2,\cdots,x_n)^\top, \bm y=(y_1,y_2,\cdots,y_m)^\top, \bm z=(z_1,z_2,\cdots,z_{m+n})^\top $.
Let $S_2=\{\bm y\in R^m| y_i>0,i=1,2,\cdots,m\}$, $S_3=R^{m+n}$ and $S=S_1\times S_2\times S_3$.
Then, a SCN form of $\Phi_1(\bm x)$  on $S$  is obtained as
\begin{eqnarray*}
g(\bm x,\bm y,\bm z)&=&\sum_{i=1}^m a_i y_i+\sum_{i=1}^m(-\ln y_i-z_i)+ \sum_{j=1}^n (-\ln x_j-z_{m+j}):\\
g_{i}(\bm x,\bm y,\bm z)&=&-\ln y_i-z_i, i=1,2,\cdots,m,\\
g_{j+m}(\bm x,\bm y,\bm z)&=&-\ln x_j-z_{m+j},j=1,2,\cdots,n,\\
h_{i}(\bm x,\bm y,\bm z)&=&y_i-\sum\limits_{j=1}^n \alpha_{ij} z_{m+j}, i=1,2,\cdots,m.
\end{eqnarray*}
\end{ex}
Hence, \eqref{d21} is equivalent to $\min\limits_{(\bm x,\bm y)}\max\limits_{\bm z} g(\bm x,\bm y,\bm z)\ \mbox{s.t.} \ (\bm x,\bm y,\bm z)\in X_c(\Phi_1)$.

\begin{ex}
Consider the support vector machine classifier via
$L_{0/1}$ Soft-Margin Loss model in Wang et al(2021). Let $(\bm x,x_0)\in R^{n+1}$  be variable.
Let the function
\begin{eqnarray}
\min\limits_{\bm x\in R^n,x_0\in R^1} \Phi_2(\bm x,x_0)=\frac{1}{2}\|\bm x \|^2+C\|(\bm 1-A\bm x-x_0\bm d)_+ \|_0,\label{d22}
\end{eqnarray}
where $A\in R^{m\times n}$, $\bm d\in R^m$ and $C>0$
are given with $\bm 1=(1,1,\cdots,1)^\top\in R^m$. Let $\bm y_1=(y_1,y_2,\cdots,y_m)^\top,\bm y_2=(y_{m+1},y_{m+2},\cdots,y_{m+m})^\top,\bm z=(z_1,z_2,\cdots,z_m)^\top,\bm y_3=(y_{2m+1},y_{2m+2},\cdots,y_{2m+m})^\top\in R^m$. Let $\bm y=(\bm y_1,\bm y_2,\bm y_3)\in R^{3m}$. Then,
a SCN form of $\Phi_2(\bm x,x_0)$ is obtained as
\begin{eqnarray*}
 g(\bm x,x_0,\bm y,\bm z)&=&\frac{1}{2}\sum\limits_{i=1}^n x_i^2+C\sum\limits_{j=1}^m [y_{m+j}^2+(y_j+y_{m+j}-1)^2-z_j+y_j^2+(y_{m+j}-1)^2-z_j]:\\
 g_j(\bm x,x_0,\bm y,\bm z)&=&(y_j+y_{m+j}-1)^2-z_j,\ j=1,2,\cdots,m,\\
 g_{m+j}(\bm x,x_0,\bm y,\bm z)&=&y_j^2+(y_{m+j}-1)^2-z_j,\ j=1,2,\cdots,m,\\
 g_{2m+j}(\bm x,x_0,\bm y,\bm z)&=& y_{m+j}^2- y_{m+j},\ j=1,2,\cdots,m,\\
 g_{3m+j}(\bm x,x_0,\bm y,\bm z)&=& y_{2m+j}^2-y_j,\ j=1,2,\cdots,m,\\
 h_{j}(\bm x,x_0,\bm y,\bm z)&=&\sum\limits_{i=1}^n a_{ij}x_i+x_0d_j+y_j-1,\ j=1,2,\cdots,m,
\end{eqnarray*}
where $S=S_1\times S_2\times S_3, S_1=R^{n+1}, S_2=\{\bm y\in R^{3n}\mid y_j\geq 0,y_{m+j}\in [0,1],j=1,2,\cdots,m\}$ and $S_3=R_+^m$.
\end{ex}%

 So, $\Phi_2(\bm x,x_0)$ is a SCN function on $S$.
Hence, \eqref{d22} is equivalent to $$\min\limits_{(\bm x,x_0,\bm y)}\max\limits_{\bm z} g(\bm x,x_0,\bm y,\bm z)\ \mbox{s.t.} \ (\bm x,x_0,\bm y,\bm z)\in X_c(\Phi_2).$$


\begin{ex}
 (In Chen et al.2010) The function in sparse optimization is
\begin{eqnarray}
\min\limits_{\bm x\in R^n} \Phi_3(\bm x)=g(\bm x)+\lambda\sum\limits_{i=1}^n \|x_i\|_0,\label{d23}
\end{eqnarray}
where $g(\bm x)$ is convex on $\bm x$. Let $\bm x,\bm z\in R^n$ and $\bm y\in R^{n}$. Then, a SCN  form of $\Phi_3(\bm x)$ on $S$ is obtained by
\begin{eqnarray*}
g(\bm x,\bm y,\bm z)&=&g(\bm x)+\lambda\sum\limits_{i=1}^n [y_{i}^2+(x_i+y_i-1)^2-z_i+x_i^2+(y_i-1)^2-z_i] :\\
g_i(\bm x,\bm y,\bm z)&=& (x_i+y_i-1)^2-z_i,\ \ i=1,2,\cdots, n,\\
g_{i+n}(\bm x,\bm y,\bm z)&=& x_i^2+(y_i-1)^2-z_i,\ \ i=1,2,\cdots, n, \\
g_{i+2n}(\bm x,\bm y,\bm z)&=&y_{i}^2-y_i,\ \ i=1,2,\cdots, n,
\end{eqnarray*}
where $S=S_1\times S_2\times S_3, S_1=R^n, S_2=\{\bm y\in R^n\mid y_i\in [0,1]\}$ and $S_3=R_+^n$.
\end{ex}%

Hence, \eqref{d23} is equivalent to $\min\limits_{(\bm x,\bm y)}\max\limits_{\bm z} g(\bm x,\bm y,\bm z)\ \mbox{s.t.} \ (\bm x,\bm y,\bm z)\in X_c(\Phi_3)$.

\begin{ex}
Let a quadratic function be
\begin{eqnarray}
\min\limits_{\bm x\in R^n}  \Phi_4(\bm x)= \bm x^\top A\bm x+\bm c^\top\bm x,\label{d24}
\end{eqnarray}
where  $\bm c=(c_1,c_2,\cdots,c_n)^\top \in R^n$ is given vectors and $A=(a_{ij})$ is a given $n\times n$ matrix. $\Phi_5(\bm x)$ is not necessarily a convex function.  Then, a SCN form of $\Phi_4(\bm x)$  is obtained by
\begin{eqnarray*}
g(\bm x,\bm z)&=& \sum\limits_{i=1}^n\sum\limits_{j=1 \atop i\not=j}^n \frac{a_{ij}sgn(a_{ij})}{2}[(x_i+sgn(a_{ij})x_j)^2-z_i-z_j]+\sum\limits_{i=1}^n [a_{ii}z_i+x_i^2-z_i+ c_ix_i]:\\
g_i(\bm x,\bm z)&=&x_i^2-z_i,\ \ \ i=1,2,\cdots,n,
\end{eqnarray*}
where $S=S_1\times S_3, S_1=R^n$ and $S_3=R_+^n$.
\end{ex}
Hence, \eqref{d24} is equivalent to $\min\limits_{(\bm x,\bm y)}\max\limits_{\bm z} g(\bm x,\bm y,\bm z)\ \mbox{s.t.} \ (\bm x,\bm y,\bm z)\in X_c(\Phi_4)$.

\begin{ex}
Let a nonconvex function in Al-Khayyaltt (1983)
\begin{eqnarray}
\min\limits_{\bm x\in R^{n+m}} \Phi_5(\bm x)=\bm c_1^\top\bm x_1 +\bm x_1^\top A\bm x_2 +\bm c_2^\top \bm x_2, \label{d25}
\end{eqnarray}
where  $\bm c_1=(c_1,c_2,\cdots,c_n)^\top\in R^n$ and $\bm c_2=(c_{n+1},c_{n+2},\cdots,c_{n+j})^\top\in R^m$  are given vectors, $A=(a_{ij})$ is a given $n\times m$ matrix, $\bm x_1=(x_1,x_2,\cdots,x_n)\top\in R^n$, $\bm x_2=(x_{n+1},x_{n+2},\cdots,x_{n+m})^\top\in R^m$ and $\bm x=(\bm x_1,\bm x_2)^\top$.
Then, a SCN form of $\Phi_5(\bm x)$  is obtained by
\begin{eqnarray*}
g(\bm x,\bm z)&=&\sum\limits_{i=1}^{n+m}(c_ix_i)+\sum\limits_{i=1}^n\sum\limits_{j=1}^m \frac{a_{ij}sgn(a_{ij})}{2}[(x_i+sgn(a_{ij})x_{n+j})^2-z_i-z_{n+j}]:\\
g_i(\bm x,\bm z)&=&x_i^2-z_i,\ \ \ i=1,2,\cdots,n,n+1,\cdots,n+m,
\end{eqnarray*}
where $S=S_1\times S_3, S_1=R^{n+m}$ and $S_3=R_+^{n+m}$.
\end{ex}%

Hence, \eqref{d25} is equivalent to $\min\limits_{(\bm x,\bm y)}\max\limits_{\bm z} g(\bm x,\bm y,\bm z)\ \mbox{s.t.} \ (\bm x,\bm y,\bm z)\in X_c(\Phi_5)$.

The above SCN functions $\Phi_i(\cdot)(i=1,2,\cdots,5)$ tell us that
$[g:g_1,g_2,\cdots,g_s;h_1,h_2,\cdots,h_r]$ is an equivalent representation of the SCN function $f$. Hence, we have $\min\limits_{\bm x\in S_1} \ f(\bm x)=\min\limits_{(\bm x,\bm y)}\max\limits_{\bm z} \{ g(\bm x,\bm y,\bm z)\mid (\bm x,\bm y,\bm z)\in X_c(f)\}$.

\section{Optimization Condition and  Exactness of (PCNO)}

In this section, let $f=[g:g_1,g_2,\cdots,g_s;h_1,h_2,\cdots,h_r]$ be defined Definition 2.1. So, $\min\limits_{\bm x\in S_1} \ f(\bm x)$ is equivalent to
\begin{eqnarray*}
\mbox{(PCNO)} &\min\limits_{(\bm x,\bm y)}\max\limits_{\bm z}& g(\bm x,\bm y,\bm z)\\
        &s.t.&    g_i(\bm x,\bm y,\bm z)\leq 0, i=1,2,\cdots,s,\\
        &&        h_j(\bm x,\bm y,\bm z)=0,j=1,2,\cdots,r,\\
        &&         (\bm x,\bm y,\bm z)\in S.
\end{eqnarray*}

\begin{theorem}
Suppose that $(\bm x^*,\bm y^*,\bm z^*)\in X(f)$ is an optimal solution to (PCNO). If there is a $(\bm x,\bm y,\bm z^*)\in X_c(f)$ and $(\bm x^*,\bm y^*,\bm z)\in X_c(f)$ such that
\begin{eqnarray}
 g_i(\bm x,\bm y,\bm z^*)< 0, g_i(\bm x^*,\bm y^*,\bm z)< 0, i=1,2,\cdots,s, \label{d26}
\end{eqnarray}
hold, then
there are $(\alpha_1,\alpha_2,\cdots,\alpha_s)^\top,(\beta_1,\beta_2,\cdots,\beta_s)^\top$,
$(\alpha_{s+1},\alpha_{s+2},\cdots,\alpha_{s+r})^\top$ and $(\beta_{s+1},\beta_{s+2},\cdots,\beta_{s+r})^\top$
such that
\begin{eqnarray}
 \nabla_{(\bm x,\bm y)} g(\bm x^*,\bm y^*,\bm z^*)+\sum\limits_{i=1}^s\alpha_i \nabla_{(\bm x,\bm y)} g_i(\bm x^*,\bm y^*,\bm z^*)+\sum\limits_{j=1}^r\alpha_{s+j} \nabla_{(\bm x,\bm y)} h_j(\bm x^*,\bm y^*,\bm z^*)=0,\label{d27} \\
 -\nabla_{\bm z} g(\bm x^*,\bm y^*,\bm z^*)+\sum\limits_{i=1}^s\beta_i \nabla_{\bm z} g_i(\bm x^*,\bm y^*,\bm z^*)+\sum\limits_{j=1}^r\beta_{s+j} \nabla_{\bm z} h_j(\bm x^*,\bm y^*,\bm z^*)=0,\label{d28} \\
          \alpha_i g_i(\bm x^*,\bm y^*,\bm z^*)=0, \alpha_i\geq 0, i=1,2,\cdots,s,\label{d29}\\
          \beta_i g_i(\bm x^*,\bm y^*,\bm z^*)=0,\beta_i\geq 0, i=1,2,\cdots,s. \label{d30}
\end{eqnarray}
\end{theorem}

\begin{ex}
Consider an  optimization problem
\begin{eqnarray*}
\mbox(P4.1)&\min& f(x)=|x|^{\frac{1}{2}}\\
&s.t.& x\in R^1.
\end{eqnarray*}
Then, a SCN optimization problem of $f(x)$ on $S$ is obtained by
\begin{eqnarray*}
\mbox(SCN4.1)&\min\limits_{x}\max\limits_{z}& g(x,y,z)=y+y^4-z+x^2-z \\
&s.t.& g_1(x,y,z)=y^4-z\leq 0,\\
&&     g_2(x,y,z)=x^2-z\leq 0,\\
&&     g_3(x,y,z)=-y\leq 0,
\end{eqnarray*}
where $S=S_1\times S_2\times S_3, S_1=R^1, S_2=R^1_+$ and $S_3=R_+^1$.

So, $x^*=0$ is an optimal solution to (P4.1). $(x^*,y^*,z^*)=(0,0,0)$ is an optimal solution to (SCN4.1).
We have $g_i(x,y,z^*)=0(i=1,2)$, $g_1(x^*,y^*,2)<0,g_2(x^*,y^*,2)<0$ and $g_3(x^*,y^*,2)=0$.
$\nabla  g(x^*,y^*,z^*)=(0,1,-2)^\top, \nabla  g_1(x^*,y^*,z^*)=(0,0,-1)^\top,
\nabla  g_2(x^*,y^*,z^*)=(0, 0, -1)^\top, \nabla  g_3(x^*,y^*,z^*)=(0, -1, 0)^\top$.
 So, when  $(\alpha_1,\alpha_2,\alpha_3)=(1,1,1)$ and  $(\beta_1,\beta_2,\beta_3)=(1,1,1)$,
 \eqref{d27},\eqref{d28},\eqref{d29} and \eqref{d30} hold.
The example shows if the condition \eqref{d26} does not hold, the conclusion of Theorem 4.1 is true.
 \end{ex}%

\begin{ex}
Consider an  optimization problem
\begin{eqnarray*}
\mbox(P4.2)&\min\limits_{x\in R^1}& f(x)=(x-1)^2+\lambda \|x\|_0,
\end{eqnarray*}
where $\lambda >0$. Then, a SCN optimization problem of $f(x)$ on $S$ is obtained by
\begin{eqnarray*}
\mbox(SCN4.2)&\min\limits_{(x,y)}\max\limits_{z}& g(x,y,z)=(x-1)^2+\lambda(y^2+(x+y-1)^2-z+x^2+(y-1)^2-z) \\
&s.t.& g_1(x,y,z)= (x+y-1)^2-z\leq 0,\\
&& g_{2}(x,y,z)= x^2+(y-1)^2-z\leq 0, \\
&& g_{3}(x,y,z)=y^2-y\leq 0, (x,y,z)\in S,
\end{eqnarray*}
where $S=S_1\times S_2\times S_3, S_1=R^1, S_2=\{y\in R^1\mid y\in [0,1]\}$ and $S_3=R_+^1$.

When $\lambda >1$, $x^*=0$ is an optimal solution to (P4.1). $(x^*,y^*,z^*)=(0,0,1)$ is an optimal solution to (SCN4.2).
We have $g_i(0.5,0.5,z^*)<0(i=1,2,3)$, $g_1(x^*,y^*,2)<0,g_2(x^*,y^*,2)<0$ and $g_3(x^*,y^*,2)=0$.
$\nabla  g(x^*,y^*,z^*)=(-2-2\lambda, -4\lambda, -2\lambda)^\top, \nabla  g_1(x^*,y^*,z^*)=(-2, -2, -1)^\top,
\nabla  g_2(x^*,y^*,z^*)=(0, -2, -1)^\top, \nabla  g_3(x^*,y^*,z^*)=(0, -1, 0)^\top$.
 So, when  $(\alpha_1,\alpha_2,\alpha_3)=(-1-\lambda,1-\lambda,0)$ and  $(\beta_1,\beta_2,\beta_3)=(2\lambda,0,0)$,
 \eqref{d26},\eqref{d27},\eqref{d28} and \eqref{d29} do not hold.

 When $0<\lambda \leq 1$, $x^*\not=0$ is an optimal solution to (P4.1). $(x^*,y^*,z^*)=(x^*,1,(x^*)^2)$ is an optimal solution to (SCN4.1).
$\nabla  g(x^*,y^*,z^*)=(2(x^*-1)+\lambda(4x^*), \lambda(2+2x^*), -2\lambda)^\top, \nabla  g_1(x^*,y^*,z^*)=(2x^*, 2x^*, -1)^\top, \nabla  g_2(x^*,y^*,z^*)=(2x^*, 0, -1)^\top, \nabla  g_3(x^*,y^*,z^*)=(0, 1, 0)^\top$.
So, when  $(\alpha_1,\alpha_2,\alpha_3)=(0,-\frac{1}{x^*}(2x^*\lambda+x^*-1)-\lambda,-(2+2x^*)\lambda)$ and  $(\beta_1,\beta_2,\beta_3)=(2\lambda,0,0)$, \eqref{d26},\eqref{d27},\eqref{d28} and \eqref{d29} do not hold.
The example shows if the condition \eqref{d26} does not hold, then the conclusion of Theorem 4.1 is not true.
 \end{ex}%

The following conclusion is clear.

\begin{theorem}
Suppose that $(\bm x^*,\bm y^*,\bm z^*)\in X(f)$. If
there are $(\alpha_1,\alpha_2,\cdots,\alpha_s)^\top,(\beta_1,\beta_2,\cdots,\beta_s)^\top$,
$(\alpha_{s+1},\alpha_{s+2},\cdots,\alpha_{s+r})^\top$ and $(\beta_{s+1},\beta_{s+2},\cdots,\beta_{s+r})^\top$
such that  \eqref{d26},\eqref{d27},\eqref{d28} and \eqref{d29} hold, then $(\bm x^*,\bm y^*,\bm z^*)$ is an optimal solution to (PCNO).
\end{theorem}

Let $\rho>0$. Two penalty functions of (PCNO) are defined by respectively
\begin{eqnarray}
F(\bm x,\bm y;\bm z,\rho)=g(\bm x,\bm y,\bm z)+\rho\sum\limits_{i=1}^s\max\{g_i(\bm x,\bm y,\bm z),0\}
   +\rho\sum\limits_{j=1}^r |h_j(\bm x,\bm y,\bm z)| \label{d31}
\end{eqnarray}
and
\begin{eqnarray}
G(\bm z;\bm x,\bm y,\rho)=-g(\bm x,\bm y,\bm z)+\rho\sum\limits_{i=1}^s\max\{g_i(\bm x,\bm y,\bm z),0\}
   +\rho\sum\limits_{j=1}^r |h_j(\bm x,\bm y,\bm z)|. \label{d32}
\end{eqnarray}
For a fixed $(\bm z,\rho)$, the corresponding optimization problem of \eqref{d31} is defined by
\begin{eqnarray*}
\mbox{PCNO}(\bm z, \rho) &\min\limits_{(\bm x,\bm y)}& F(\bm x,\bm y;\bm z,\rho)\\
        &s.t.&            (\bm x,\bm y,\bm z)\in S.
\end{eqnarray*}
For a fixed $(\bm x,\bm y,\rho)$, the corresponding optimization problem of \eqref{d32} is defined by
\begin{eqnarray*}
\mbox{PCNO}(\bm x,\bm y, \rho) &\min\limits_{\bm z}& G(\bm z;\bm x,\bm y,\rho)\\
        &s.t.&            (\bm x,\bm y,\bm z)\in S.
\end{eqnarray*}

\begin{definition}
Suppose that $(\bm x^*,\bm y^*,\bm z^*)\in X(f)$ is an optimal solution to (PCNO).
If there is a $\rho'>0$ such that $(\bm x^*,\bm y^*)$ is an optimal solution to (PCNO)$(\bm z^*, \rho)$ for $\rho>\rho'$ and $\bm z^*$ is an optimal solution to (PCNO)$(\bm x^*,\bm y^*, \rho)$ for $\rho>\rho'$, then
\eqref{d31} and  \eqref{d32} are is exact.
\end{definition}

\begin{theorem}
Suppose that $(\bm x^*,\bm y^*,\bm z^*)\in X(f)$ is an optimal solution to (PCNO). If
there are $(\alpha_1,\alpha_2,\cdots,\alpha_s)^\top,(\beta_1,\beta_2,\cdots,\beta_s)^\top$,
$(\alpha_{s+1},\alpha_{s+2},\cdots,\alpha_{s+r})^\top$ and $(\beta_{s+1},\beta_{s+2},\cdots,\beta_{s+r})^\top$
such that  \eqref{d26}),\eqref{d27},\eqref{d28} and \eqref{d29} hold, then \eqref{d31} and  \eqref{d32} are exact for $\rho>\rho^*=\max\{\alpha_i,\beta_i| i=1,2,\cdots,s+r\}$.
\end{theorem}
{\bf Proof.} Since $g,g_i(i=1,2,\cdots,s)$ are convex on $(\bm x,\bm y)$ and $h_j(j=1,2,\cdots,r)$ are linear on $(\bm x,\bm y)$,
we have
\begin{eqnarray}
g(\bm x,\bm y,\bm z^*)-g(\bm x^*,\bm y^*,\bm z^*)
\geq \nabla_{(\bm x,\bm y)} g(\bm x^*,\bm y^*,\bm z^*)[(\bm x,\bm y,\bm z^*)-\bm (x^*,\bm y^*,\bm z^*)],\label{d33}\\
g_i(\bm x,\bm y,\bm z^*)-g_i(\bm x^*,\bm y^*,\bm z^*)\geq \nabla_{(\bm x,\bm y)} g_i(\bm x^*,\bm y^*,\bm z^*)[(\bm x,\bm y,\bm z^*)-\bm (x^*,\bm y^*,\bm z^*)],i=1,2,\cdots,s,\label{d34}\\
h_j(\bm x,\bm y,\bm z^*)-h_j(\bm x^*,\bm y^*,\bm z^*)= \nabla_{(\bm x,\bm y)} h_j(\bm x^*,\bm y^*,\bm z^*)[(\bm x,\bm y,\bm z^*)-\bm (x^*,\bm y^*,\bm z^*)],i=1,2,\cdots,r.\label{d35}
\end{eqnarray}
For $(\bm x,\bm y,\bm z^*)\in X_c(f)$ and $\rho>\rho^*=\max\{\alpha_i,\beta_i| i=1,2,\cdots,s+r\}$, by \eqref{d33},\eqref{d34} and \eqref{d35}, we have
\begin{eqnarray*}
F(\bm x,\bm y;\bm z^*,\rho)-F(\bm x^*,\bm y^*,\bm z^*,\rho)&=&
g(\bm x,\bm y,\bm z^*)-g(\bm x^*,\bm y^*,\bm z^*)\\
&&+\rho\sum\limits_{i=1}^s(\max\{g_i(\bm x,\bm y,\bm z^*),0\}-\max\{g_i(\bm x^*,\bm y^*,\bm z^*),0\})\\
&&+\rho\sum\limits_{j=1}^r(|h_j(\bm x,\bm y,\bm z^*)|-|h_j(\bm x^*,\bm y^*,\bm z^*)|)\\
&\geq& \nabla_{(\bm x,\bm y)} g(\bm x^*,\bm y^*,\bm z^*)[(\bm x,\bm
y,\bm z^*)-\bm (x^*,\bm y^*,\bm z^*)],\\
&&+\sum\limits_{i=1}^s \alpha_i\max\{g_i(\bm x,\bm y,\bm z^*),0\}
+\sum\limits_{j=1}^r\alpha_{s+j}|h_j(\bm x,\bm y,\bm z^*)|\\
&\geq& \nabla_{(\bm x,\bm y)} g(\bm x^*,\bm y^*,\bm z^*)[(\bm x,\bm
y,\bm z^*)-\bm (x^*,\bm y^*,\bm z^*)],\\
&&+\sum\limits_{i=1}^s \alpha_i(g_i(\bm x,\bm y,\bm z^*)-g_i(\bm x^*,\bm y^*,\bm z^*))\\
&& +\sum\limits_{j=1}^r\alpha_{s+j}(h_j(\bm x,\bm y,\bm z^*)-h_j(\bm x^*,\bm y^*,\bm z^*))
\end{eqnarray*}
since $g_i(\bm x^*,\bm y^*,\bm z^*)\leq 0,\alpha_ig_i(\bm x^*,\bm y^*,\bm z^*)=0,i=1,2,\cdots,i,$ $h_j(\bm x^*,\bm y^*,\bm z^*)=0, j=1,2,\cdots,r$.
By \eqref{d26},\eqref{d27},\eqref{d28} and \eqref{d29}, we have
\begin{eqnarray*}
F(\bm x,\bm y;\bm z^*,\rho)-F(\bm x^*,\bm y^*,\bm z^*,\rho)&\geq& \nabla_{(\bm x,\bm y)} g(\bm x^*,\bm y^*,\bm z^*)[(\bm x,\bm y,\bm z^*)-\bm (x^*,\bm y^*,\bm z^*)],\\
&&+\sum\limits_{i=1}^s \alpha_i\nabla_{(\bm x,\bm y)} g_i(\bm x^*,\bm y^*,\bm z^*)[(\bm x,\bm y,\bm z^*)-\bm (x^*,\bm y^*,\bm z^*)]\\
&&+\sum\limits_{j=1}^r\alpha_{s+j}\nabla_{(\bm x,\bm y)} h_j(\bm x^*,\bm y^*,\bm z^*)[(\bm x,\bm y,\bm z^*)-\bm (x^*,\bm y^*,\bm z^*)]\\
&\geq& 0.
\end{eqnarray*}
Similarly, $G(\bm z;\bm x^*,\bm y^*,\rho)\geq G(\bm z^*;\bm x^*,\bm y^*,\rho)$ for $\rho>\rho^*$.

\begin{theorem}
Suppose that $(\bm x^*,\bm y^*,\bm z^*)\in X(f)$ is an optimal solution to (PCNO). If \eqref{d31} and  \eqref{d32} are exact, then
there are $(\alpha_1,\alpha_2,\cdots,\alpha_s)^\top,(\beta_1,\beta_2,\cdots,\beta_s)^\top$,
$(\alpha_{s+1},\alpha_{s+2},\cdots,\alpha_{s+r})^\top$ and $(\beta_{s+1},\beta_{s+2},\cdots,\beta_{s+r})^\top$
such that  \eqref{d26},\eqref{d27},\eqref{d28} and \eqref{d29} hold.
\end{theorem}
{\bf Proof.} By \cite{Clarke}, the conclusion is true.

The perturbation set of $X_c(f)$ is defined by
\begin{eqnarray}
  X_c(f,\bm \eta,\bm \tau)=\{(\bm x,\bm y,\bm z)\in S&\mid & g_i(\bm x,\bm y,\bm z)\leq \eta_i,i=1,2,\cdots,s;\nonumber\\
   &&h_j(\bm x,\bm y,\bm z)=\tau_j, j=1,2,\cdots,r\}, \label{d351}
\end{eqnarray}
where $\bm \eta=(\eta_1,\eta_2,\cdots, \eta_s)^\top\in R^s$ and $\bm \tau=(\tau_1,\tau_2,\cdots,\tau_r)^\top\in R^r$.

A perturbed problem of (PCNO) is defined by
\begin{eqnarray*}
\mbox{PCNO}(\bm \eta,\bm \tau) &\min\limits_{(\bm x,\bm y)}&\max\limits_{\bm z} g(\bm x,\bm y,\bm z)\\
        &s.t.&          (\bm x,\bm y,\bm z)\in  X_c(f,\bm \eta,\bm \tau).
\end{eqnarray*}

\begin{definition}\label{de2.1}
Let $(\bm x^*,\bm y^*,\bm z^*)$ be an optimal solution to (PCNO) and
 $(\bm x^*_{(\bm \eta,\bm \tau)},\bm y^*_{(\bm \eta,\bm \tau)},\bm z^*_{(\bm \eta,\bm \tau)})$ be an optimal solution to PCNO$(\bm \eta,\bm \tau)$
  for any $(\bm \eta,\bm \tau)\in R^s\times R^r$.
If there is a $\rho'$ such that
\begin{eqnarray}
  |g(\bm x^*,\bm y^*,\bm z^*)-g(\bm x^*_{(\bm \eta,\bm \tau)},\bm y^*_{(\bm \eta,\bm \tau)},\bm z^*_{(\bm \eta,\bm \tau)})|\leq \rho |(\bm \eta,\bm \tau)|, \ \ \ \ \forall  \rho>\rho', \label{d352}
  \end{eqnarray}
where $|(\bm \eta,\bm \tau)|=\sum\limits_{i=1}^{s}\max\{\eta_i,0\}+\sum\limits_{j=1}^{r}|\tau_j|$, then $F(\bm x,\bm y;\bm z,\rho)$ and $G(\bm z;\bm x,\bm y,\rho)$ are  stable.
\end{definition}

\begin{theorem}\label{th2.6}
Let $(\bm x^*,\bm y^*,\bm z^*)$ be an optimal solution to (PCNO). Then,
$F(\bm x,\bm y;\bm z,\rho)$ and $G(\bm z;\bm x,\bm y,\rho)$ are stable  if  and only if
$F(\bm x,\bm y;\bm z,\rho)$ and $G(\bm z;\bm x,\bm y,\rho)$ are exact.
\end{theorem}

{\bf Proof.} That $F(\bm x,\bm y;\bm z,\rho)$ and $G(\bm z;\bm x,\bm y,\rho)$ are exact  when  $F(\bm x,\bm y;\bm z,\rho)$ and $G(\bm z;\bm x,\bm y,\rho)$ are stable is proved first. According to the Definition 4.2, for any $(\bm \eta,\bm \tau)$, we obtain that there is a $\rho'$ satisfying that \eqref{d352} holds. Suppose that $F(\bm x,\bm y;\bm z,\rho)$ and $G(\bm z;\bm x,\bm y,\rho)$ are not exact.
 Then, there always exist some $\rho>\rho'$ and $(\bm x',\bm y',\bm z')$ such that
$$F(\bm x',\bm y';\bm z^*,\rho)<F(\bm x^*,\bm y^*;\bm z^*,\rho)=g(\bm x^*,\bm y^*;\bm z^*),$$
$$G(\bm z',;\bm x^*,\bm y^*,\rho)<G(\bm z^*;\bm x^*,\bm y^*,\rho)=-g(\bm x^*,\bm y^*;\bm z^*).$$
Thus,
\begin{eqnarray*}
g(\bm x',\bm y',\bm z^*)+\rho\sum\limits_{i=1}^s\max\{g_i(\bm x',\bm y',\bm z^*),0\}
+\rho\sum\limits_{j=1}^r |h_j(\bm x',\bm y',\bm z^*)|<g(\bm x^*,\bm y^*;\bm z^*)\\
-g(\bm x^*,\bm y^*,\bm z')+\rho\sum\limits_{i=1}^s\max\{g_i(\bm x^*,\bm y^*,\bm z'),0\}
   +\rho\sum\limits_{j=1}^r |h_j(\bm x^*,\bm y^*,\bm z')|<-g(\bm x^*,\bm y^*;\bm z^*).
\end{eqnarray*}
If $(\bm x',\bm y',\bm z^*)\in X_c(f)$ and $(\bm x^*,\bm y^*,\bm z')\in X_c(f)$, then
\begin{eqnarray*}
g(\bm x',\bm y',\bm z^*)<g(\bm x^*,\bm y^*;\bm z^*)<g(\bm x^*,\bm y^*,\bm z').
\end{eqnarray*}
This implies that $(\bm x^*,\bm y^*,\bm z^*)$ is not an optimal solution to (PCNO). A contradiction occurs. Hence,  $(\bm x',\bm y',\bm z^*)\in X_c(f)$ and $(\bm x^*,\bm y^*,\bm z')\in X_c(f)$ do not hold,
 and $\rho\sum\limits_{i=1}^s\max\{g_i(\bm x',\bm y',\bm z^*),0\}
+\rho\sum\limits_{j=1}^r |h_j(\bm x',\bm y',\bm z^*)|>0$ or $\rho\sum\limits_{i=1}^s\max\{g_i(\bm x^*,\bm y^*,\bm z'),0\}   +\rho\sum\limits_{j=1}^r |h_j(\bm x^*,\bm y^*,\bm z')|)>0$. Well, let
$$\rho\sum\limits_{i=1}^s\max\{g_i(\bm x',\bm y',\bm z^*),0\}+\rho\sum\limits_{j=1}^r |h_j(\bm x',\bm y',\bm z^*)|>0$$ be established. Then,
let $\bm \eta'=(\eta_1',\eta_2',\cdots, \eta_s')^\top\in R^s$ and $\bm \tau'=(\tau_1',\tau_2',\cdots,\tau_r')^\top\in R^r$ with $\eta_i'=\max\{g_i(\bm x',\bm y',\bm z^*),0\}$ $(i=1,2,\cdots,s)$ and
$\tau_j'=|h_j(\bm x',\bm y',\bm z^*)$ $(j=1,2,\cdots,r)|$
, and $(\bm x^*_{(\bm \eta',\bm \tau')},\bm y^*_{(\bm \eta',\bm \tau')},\bm z^*_{(\bm \eta',\bm \tau')})$ be an optimal solution to PCNO$(\bm \eta',\bm \tau')$.
Then,
$$g(\bm x^*_{(\bm \eta',\bm \tau')},\bm y^*_{(\bm \eta',\bm \tau')},\bm z^*_{(\bm \eta',\bm \tau')})\leq g(\bm x',\bm y',\bm z^*).$$
Therefore,
\begin{eqnarray*}
g(\bm x^*_{(\bm \eta',\bm \tau')},\bm y^*_{(\bm \eta',\bm \tau')},\bm z^*_{(\bm \eta',\bm \tau')})+\rho\sum\limits_{i=1}^s \eta_i'+\rho\sum\limits_{j=1}^r \tau_j'
&\leq& g(\bm x',\bm y',\bm z^*)+\rho\sum\limits_{i=1}^s \eta_i'+\rho\sum\limits_{j=1}^r \tau_j' \\
&=&F(\bm x',\bm  y';\bm z^*,\rho)<g(\bm x^*,\bm y^*;\bm z^*),
\end{eqnarray*}
which shows that
$$g(\bm x^*,\bm y^*;\bm z^*)-g(\bm x^*_{(\bm \eta',\bm \tau')},\bm y^*_{(\bm \eta',\bm \tau')},\bm z^*_{(\bm \eta',\bm \tau')})>\rho|(\bm \eta',\bm \tau')|,$$
where  $|(\bm \eta',\bm \tau')|=\sum\limits_{i=1}^s \eta_i'+\sum\limits_{j=1}^r \tau_j'$. This inequality contradicts  \eqref{d352}. Hence,  $F(\bm x,\bm y;\bm z,\rho)$ and $G(\bm z;\bm x,\bm y,\rho)$ are not stable,
which yields a contradiction with the assumption and proves that $F(\bm x,\bm y;\bm z,\rho)$ and $G(\bm z;\bm x,\bm y,\rho)$ are  exact.

 Next, that $F(\bm x,\bm y;\bm z,\rho)$ and $G(\bm z;\bm x,\bm y,\rho)$ are stable when $F(\bm x,\bm y;\bm z,\rho)$ and $G(\bm z;\bm x,\bm y,\rho)$ are exact  is proved. According to the definition of Definition 4.1,
there is a $\rho'>0$ such that $(\bm x^*,\bm y^*)$ is an optimal solution to (PCNO)$(\bm z^*, \rho)$ for $\rho>\rho'$ and $\bm z^*$ is an optimal solution to (PCNO)$(\bm x^*,\bm y^*, \rho)$ for $\rho>\rho'$, i.e.
\begin{eqnarray}
F({\bm x}^*, {\bm y}^*; \bm z^*,\rho)\leq F(\bm x, {\bm y};\bm z^*,\rho),  \forall (\bm x,\bm y,\bm z^*)\in S,\label{d36}\\
G({\bm z}^*; {\bm x}^*,\bm y^*,\rho)\leq G({\bm z};{\bm x}^*,\bm y^*,\rho),  \forall (\bm x^*,\bm y^*,\bm z)\in S.\label{d37}
\end{eqnarray}
Let $(\bm x^*_{(\bm \eta,\bm \tau)},\bm y^*_{(\bm \eta,\bm \tau)},\bm z^*_{(\bm \eta,\bm \tau)})$ be an optimal solution to PCNO$(\bm \eta,\bm \tau)$  for any $(\bm \eta,\bm \tau)\in R^s\times R^r$. By \eqref{d36} and \eqref{d37}, we have
\begin{eqnarray}
g({\bm x}^*, {\bm y}^*, \bm z^*)=F({\bm x}^*, {\bm y}^*; \bm z^*,\rho)\leq F(\bm x^*_{(\bm \eta,\bm \tau)},\bm y^*_{(\bm \eta,\bm \tau)};\bm z^*_{(\bm \eta,\bm \tau)},\rho), \label{d38} \\
-g({\bm x}^*, {\bm y}^*, \bm z^*)=G({\bm z}^*; {\bm x}^*,\bm y^*,\rho)\leq G(\bm z^*_{(\bm \eta,\bm \tau)};\bm x^*_{(\bm \eta,\bm \tau)},\bm y^*_{(\bm \eta,\bm \tau)},\rho).\label{d39}
\end{eqnarray}
From \eqref{d38} and \eqref{d39}, thus
$$  |g({\bm x}^*, {\bm y}^*, \bm z^*)-g(\bm x^*_{(\bm \eta,\bm \tau)},\bm y^*_{(\bm \eta,\bm \tau)},\bm z^*_{(\bm \eta,\bm \tau)})\leq \rho |(\bm \eta,\bm \tau)|, \ \ \ \ \forall  \rho>\rho'.$$
It follows from the definition that \eqref{d31} and \eqref{d32} are stable.\\

Theorem 4.5 shows that a approximate optimal solution to (PCNO) may be obtained by by solving the following penalty function problem {PCNO}$(\bm z, \rho)$ and {PCNO}$(\bm x,\bm y, \rho)$, when \eqref{d31} and \eqref{d32} are stable or exact.

For a fixed $(\bm x,\bm y,\rho)$, a penalty function is defined by
\begin{eqnarray}
G_2(\bm z;\bm x,\bm y,\rho)=-g(\bm x,\bm y,\bm z)+\rho\sum\limits_{i=1}^s\max\{g_i(\bm x,\bm y,\bm z),0\}^2
   +\rho\sum\limits_{j=1}^r h_j(\bm x,\bm y,\bm z)^2. \label{d40}
\end{eqnarray}
Since $g_i(i=1,2,\cdots,s)$ is convex on $\bm z$ and  $h_j(j=1,2,\cdots,r)$ is linear, $G_2(\bm z;\bm x,\bm y,\rho)$ is convex on $\bm z$.
For a fixed $(\bm x,\bm y,\rho)$, the corresponding optimization problem of \eqref{d40} is defined by
\begin{eqnarray*}
\mbox{PCNO}_2(\bm x,\bm y, \rho) &\min\limits_{\bm z}& G_2(\bm z;\bm x,\bm y,\rho)\\
        &s.t.&            (\bm x,\bm y,\bm z)\in S.
\end{eqnarray*}
The following conclusion is clear.

\begin{theorem}\label{th2.6}
For a given $(\bm x,\bm y,\bm z^*)\in X_c(f)$, if $\bm z^*$ is an optimal solution to PCNO$_2(\bm x,\bm y, \rho)$,then $\nabla_{\bm z} G_2(\bm z^*;\bm x,\bm y,\rho)=0$ and $\bm z^*\in \arg\max\limits{\bm z}\{ g((\bm x,\bm y,\bm z)\mid (\bm x,\bm y,\bm z)\in X_c(f)\}$.
\end{theorem}

The SCN forms of the SCN functions of many examples in Section 3 show that  $g(\bm x,\bm y,\bm z)$ and $g_i(\bm x,\bm y,\bm z)(i=1,2,\cdots,r)$ is convex on $(\bm x,\bm y)$ and linear on $\bm z$. So, $\nabla_{\bm z} G_2(\bm z^*;\bm x,\bm y,\rho)$ is liner on $\bm z$.

\section{Exact Penalty Function Algorithm of (PCNO)}

In this section, the  algorithm  to approximate optimal solution to (PCNO) is studied by solving the following penalty function problem {PCNO}$(\bm z, \rho)$ and {PCNO}$(\bm x,\bm y, \rho)$ for a fixed $(\bm x,\bm y,\bm z)$.
\begin{definition} \label{de2.1}
Let $(\bm x,\bm y,\bm z)\in S$ and $\epsilon>0$. Then $(\bm x,\bm y,\bm z)$ is called an $\epsilon$-feasible solution to (PCNO), if
\begin{eqnarray}
  X_c(f,\epsilon)=\{(\bm x,\bm y,\bm z)\in S&\mid & g_i(\bm x,\bm y,\bm z)\leq \epsilon,i=1,2,\cdots,s;\nonumber\\
   &&|h_j(\bm x,\bm y,\bm z)|\leq \epsilon, j=1,2,\cdots,r\}, \label{d40}
\end{eqnarray}
\end{definition}

Based on above results, a generic algorithm to compute approximately optimal solution to (PCNO) is presented. The algorithm is called Approximate Penalty Function Algorithm for (PCNO) (APFA for short). \vspace{5mm}
\\

{\bf APFA Algorithm:}
\begin{description}
\item[Step 1:] Choose $\rho_1\geq 1, N>1$, $(\bm x^1,\bm y^1,\bm z^1)$ and $k=1$.
\item[Step 2:] Let $(\bm x^k,\bm y^k,\bm z^k)$ be obtained.
For a fixed $(\bm z^k,\rho_k)$, solve $(\bm x^{k+1},\bm y^{k+1})$ to be an optimal solution to
\begin{eqnarray*}
\min\limits_{(\bm x,\bm y,\bm z^k)\in S}~~ F(\bm x,\bm y;\bm z^k,\rho_k).
\end{eqnarray*}
Then, for a fixed $(\bm x^{k+1},\bm y^{k+1},\rho_k)$,  solve $\bm z^{k+1}$ to be an optimal solution to
\begin{eqnarray*}
\min\limits_{(\bm x^{k+1},\bm y^{k+1},\bm z)\in S}~~ G(\bm x^{k+1},\bm y^{k+1},\rho_k).
\end{eqnarray*}
\item[Step 3:] If $\|(\bm x^k,\bm y^k,\bm z^k)-(\bm x^{k+1},\bm y^{k+1},\bm z^{k+1})\|\leq \epsilon$ and $(\bm x^{k+1},\bm y^{k+1},\bm z^{k+1})\in X_c(f,\epsilon)$,
stop and $(\bm x^{k+1},\bm y^{k+1},\bm z^{k+1})$ is $\epsilon$-feasible solution to (PCNO).
          Otherwise, let $\rho_{k+1}=N\rho_k, k=:k+1$ and go to Step 2.
\end{description}

 The convergence of the APFA algorithm
is proved in the following theorem. Let
$$S(L, g)=\{(\bm x,\bm y,\bm z)\mid g(\bm x,\bm y,\bm z)\leq L, \; k=1, 2, \cdots\}, $$
which is called a Q-level set. We say that $S(L, f)$ is bounded if $S(L, f)$ is bounded for any given $L>0$ .

\begin{theorem}
\label{th3.1}
Let $\{(\bm x^k,\bm y^k,\bm z^k)\}$ be the
sequence generated by the APFA algorithm and $S$ be compact.

(i) If $\{(\bm x^k,\bm y^k,\bm z^k)\}(k=1, 2, \cdots, \bar{k})$ is a finite sequence (i.e., the APFA algorithm stops at the $\bar{k}$-th iteration) for $k>1$, then $(\bm x_1^{\bar{k}},\bm x_2^{\bar{k}})$ is $\epsilon$-feasible solution to (PCNO).
Furthermore, $(\bm x^k,\bm y^k,\bm z^k)=(\bm x^{k+1},\bm y^{k+1},\bm z^{k+1})$ and $(\bm x^{k+1},\bm y^{k+1},\bm z^{k+1})\in X_c(f)$, then $(\bm x_1^{\bar{k}},\bm x_2^{\bar{k}})$ is an optimal solution to (PCNO).

(ii) Let $\{(\bm x^k,\bm y^k,\bm z^k)\}$ be an infinite sequence,
sequence $\{F(\bm x^k,\bm y^k;\bm z^k,\rho_k)\}$ and $\{G(\bm z^k;\bm x^k,\bm y^k,\rho_k)\}$ be bounded and the Q-level set $S(L,g)$ be bounded. Then $\{(\bm x^k,\bm y^k;\bm z^k)\}$ is bounded and any limit point $(\bm x^*,\bm y^*;\bm z^*)$ of it is an optimal solution to (PCNO).
\end{theorem}

{\bf Proof.}  (i) The conclusion is clear.

(iii) By the APFA algorithm, since $\{F(\bm x^k,\bm y^k;\bm z^k,\rho_k)\}$ and $\{G(\bm z^k;\bm x^k,\bm y^k,\rho_k)\}$ are bounded as $k\to +\infty$,  there must be some $L>0$ and $k'>1$ such
that\begin{eqnarray*}
          L&>&F(\bm x^k,\bm y^k;\bm z^k,\rho_k)\\
          &\geq& g(\bm x^k,\bm y^k,\bm z^k)+\rho_k P(\bm x^k,\bm y^k,\bm z^k)\\
          &\geq& g(\bm x^k,\bm y^k,\bm z^k), ~~\forall k,
\end{eqnarray*}
and \begin{eqnarray*}
          L&>&G(\bm z^k;\bm x^k,\bm y^k,\rho_k)\\
          &\geq& -g(\bm x^k,\bm y^k,\bm z^k)+\rho_k P(\bm x^k,\bm y^k,\bm z^k)\\
          &\geq& -g(\bm x^k,\bm y^k,\bm z^k), ~~\forall k,
\end{eqnarray*}
where $P(\bm x^k,\bm y^k,\bm z^k)=\sum\limits_{i=1}^s\max\{g_i(\bm x^k,\bm y^k,\bm z^k),0\}
   +\sum\limits_{j=1}^r |h_j(\bm x^k,\bm y^k;\bm z^k)|$.
The conclusion is that $\{(\bm x^k,\bm y^k,\bm z^k)\}$ is bounded because the Q-level set $S(L,g)$ is bounded. We have
$\frac{2L}{\rho_k}>P(\bm x^k,\bm y^k,\bm z^k)\to 0$ as  $k\to +\infty$.
Without loss of generality, suppose $(\bm x^k,\bm y^k,\bm z^k)\to (\bm x^*,\bm y^*,\bm z^*)$. So, $(\bm x^*,\bm y^*,\bm z^*)$ is a feasible solution to (PCNO).

By the Step 2 of APFA Algorithm, we have
$$F(\bm x^{k+1},\bm y^{k+1};\bm z^k,\rho_k)\leq F(\bm x,\bm y;\bm z^k,\rho_k), \forall (\bm x,\bm y,\bm z^k)\in S,$$
$$F(\bm z^{k+1};\bm x^{k+1},\bm y^{k+1},\rho_k)\leq G(\bm z;\bm x^{k+1},\bm y^{k+1},\rho_k), \forall (\bm x^{k+1},\bm y^{k+1},\bm z)\in S,$$
where $k>k'$. For $(\bm x,\bm y;\bm z^k,\rho_k)\in X_c(f)$ and $(\bm z;\bm x^{k+1},\bm y^{k+1})\in X_c(f)$ let $k\to +\infty$, the above inequations are
$$F(\bm x^{*},\bm y^{*};\bm z^*,\rho_k)\leq F(\bm x,\bm y;\bm z^*,\rho_k), \forall (\bm x,\bm y,\bm z^*)\in X_c(f),$$
$$F(\bm z^{*};\bm x^{*},\bm y^*,\rho_k)\leq G(\bm z;\bm x^{*},\bm y^{*},\rho_k), \forall (\bm x^{*},\bm y^{*},\bm z)\in X_c(f),$$
i.e.
$$g(\bm x^{*},\bm y^{*},\bm z^*)\leq g(\bm x,\bm y;\bm z^*), \forall (\bm x,\bm y;\bm z^*)\in X_c(f),$$
$$-g(\bm x^{*},\bm y^*,\bm z^{*})\leq g(\bm x^{*},\bm y^{*},\bm z), \forall (\bm x^{*},\bm y^{*},\bm z)\in X_c(f).$$
Hence, $(\bm x^{*},\bm y^{*},\bm z^*)$ is an optimal solution to (PCNO).
\hfill\fbox{}

Since \eqref{d31} and  \eqref{d32} are nonsmooth, define approximate penalty function of \eqref{d31} and  \eqref{d32},
\begin{eqnarray*}
F_\theta(\bm x,\bm y;\bm z,\rho)=g(\bm x,\bm y,\bm z)+\rho \sum\limits_{i=1}^s g_i^+(\bm x,\bm y,\bm z)^\theta+\rho \sum\limits_{j=1}^r |h_j(\bm x,\bm y,\bm z)|^\theta,\nonumber\\
G_\theta(\bm z;\bm x,\bm y,\rho)=-g(\bm x,\bm y,\bm z)+\rho \sum\limits_{i=1}^s g_i^+(\bm x,\bm y,\bm z)^\theta+\rho \sum\limits_{j=1}^r |h_j(\bm x,\bm y,\bm z)|^\theta,
\end{eqnarray*}
where $\rho>0, \theta>1$ is a penalty parameter and $g_i^+(\bm x,\bm y,\bm z)=\max\{g_i(\bm x,\bm y,\bm z),0\}$. $\theta$ closes to 1. In APFA algorithm, \eqref{d31} and  \eqref{d32} are replaced with $F_\theta(\bm x,\bm y;\bm z,\rho)$ and $G_\theta(\bm x,\bm y;\bm z,\rho)$, then we solve approximate solution to (PCNO) with Matlab.

\begin{ex}
A nonconvex nonsmooth optimization problem is (Problem 5 in \cite{Bagirov})
 \begin{eqnarray*}
\mbox{(P5.1)} &\min& f_n(\bm x)=n\max\{|x_i|: i=1,2,\cdots,n\}-\sum\limits_{i=1}^n |x_i| \\
&s.t.& \bm x\in R^n.
\end{eqnarray*}
An optimal solution to (P5.1) is $\bm x^*=(\pm\alpha,\pm\alpha,\cdots,\pm\alpha)^\top$ in \cite{Bagirov} with $f(\bm x^*)=0$ for $\alpha \in R^1$ in \cite{Bagirov}. Let $\bm x,\bm z\in R^n,\bm y\in R^{n+1}$. A SCN form of $f_n$ is defined by
 \begin{eqnarray*}
g(\bm x,\bm y,\bm z)&=&n y_{n+1}-\sum\limits_{i=1}^n (y_i-y_i^2+z_{i}-x_i^2+z_{i}): \\
g_i(\bm x,\bm y,\bm z)&=& y_i^2-z_{i}\leq 0,\ i=1,2,\cdots,n,\\
g_{n+i}(\bm x,\bm y,\bm z)&=&x_i^2-z_{i}\leq 0,\ i=1,2,\cdots,n,\\
g_{2n+i}(\bm x,\bm y,\bm z)&=&y_i-y_{n+1}\leq 0, \ i=1,2,\cdots,n,
\end{eqnarray*}
where $S=R^n\times R^n_+\times R^n_+$.
In  ALPF Algorithm, the starting parameters $\epsilon=10^{-6},\rho_{1}=10,N=100,\theta=1.01$ and
$(\bm x^1,\bm y^1)=(1, 2, 3,\cdots,2n+1)$ and $\bm z^1=(-1000,-2000,\cdots,-n*1000)$ are taken.  For $n=5$, at the 1th step, an approximate solution
$\bm x^2=(0.002684,0.002684,0.002684,$ $ 0.002684,0.002684)$ to (P5.1) is obtained. In  ALPF Algorithm, the starting parameters $\epsilon=10^{-6},\rho_{1}=10,N=100,\theta=1.01$ and
$(\bm x^1,\bm y^1)=(1, 2, 3,\cdots,2n+1)$ and $\bm z^1=(1000,2000,\cdots,n*1000)$ are taken.
  For $n=10$, at the 2th step, an approximate solution $\bm x^3=(20.000000,20.000000,$ $20.000000, 20.000000,20.000000,20.000000,20.000000,20.000000, 20.000000,20.000000)$ to (P5.1)
is obtained.
\end{ex}

Hence, the above examples  illustrate that it is efficient to solve an approximate optimal solution to (PCNO) by using the ALPF algorithm with Maltlab, to avoids using subdifferentiation.

\section{Conclusion}

This paper solves a difficulty relating to unconstrained nonconvex nonsmooth optimization problems. The unconstrained, nonconvex and nonsmooth optimization problem is transformed into a convex constrained optimization problem with a convex-concave and smooth objective function.  A new concept - SCN function - is proposed, which covers many nonconvex or nonsmooth functions, even discontinuous nonconvex functions. Some of SCN function's operational properties are proved. The SCN forms of many important functions are given.A minmax optimization problem of a SCN function is defined. The equivalence of optimality condition, exactness and stability of this minmax problem are proved.

This paper provides a feasible idea for solving nonconvex or nonsmooth optimization problems with SCN functions, which shows its potential importance in solving these problems in many application fields. The major advantage of SCN function is that  the objective function of the transformed minmax problem is  convex concave  and the constraint set is  convex. Because there are many good algorithms to solve such minmax problems, the SCN function technique makes itself an effective method to solve nonconvex and nonsmooth functions.

There are at least three directions worthy of further study in terms of SCN optimization problems:

(1) decomposable Newton algorithm or decomposable SQP algorithm,

(2) Lagrangian multiplier alternating algorithm,

(3) some special forms of SCN optimization problems, for example when $f$ is a SCN function with $f=[g,g_1,g_2,\cdots,g_r]$, where  $g,g_1,g_2,\cdots,g_r$ are quadratic and linear.

{\bf Acknowledgments}

This work is supported by the National Natural Science Foundation of China(No.11871434) and the Natural Science Foundation of Zhejiang Province(No.LY18A010031).

\small
\bibliographystyle{siamplain}

\end{document}